\documentclass{amsart}

\usepackage{amsfonts}
\usepackage{amsmath}
\usepackage{amsthm}
\usepackage{amssymb}
\usepackage[english]{babel}
\usepackage{graphics}
\usepackage{latexsym}
\usepackage{longtable} 
\usepackage{mathrsfs} 
\usepackage{graphicx}
\usepackage{wrapfig} 
\usepackage[pdftex,colorlinks=true,linkcolor=blue,citecolor=magenta]{hyperref}
\usepackage{enumitem}
\usepackage{esint}

\newtheorem{theorem}{Theorem}
\newtheorem{corollary}[theorem]{Corollary}
\newtheorem{lemma}[theorem]{Lemma}
\newtheorem{proposition}[theorem]{Proposition}

\newtheorem{claim}[theorem]{Claim}
\newtheorem{example}[theorem]{Example}

\theoremstyle{definition}
\newtheorem{definition}[theorem]{Definition}
\newtheorem{remark}[theorem]{Remark}

\newcommand{\D}{\mathrm{D}}

\newcommand{\F}{\mathrm{F}}
\newcommand{\G}{\mathrm{G}}

\newcommand{\R}{\mathbb{R}}

\newcommand{\N}{\mathbb{N}}

\newcommand{\mS}{\mathbb{S}}
\newcommand{\mB}{\mathbb{B}}

\newcommand{\X}{{\bf X}}
\newcommand{\Y}{{\bf Y}}

\renewcommand{\P}{\mathrm{P}}
\newcommand{\Q}{\mathrm{Q}}

\newcommand{\noi}{\noindent}
\newcommand{\ms}{\medskip}

\newcommand{\al}{\alpha}
\newcommand{\be}{\beta}
\newcommand{\ga}{\gamma}

\newcommand{\de}{\delta}
\newcommand{\De}{\Delta}
\newcommand{\e}{\varepsilon}
\newcommand{\si}{\sigma}

\newcommand{\la}{\lambda}

\newcommand{\Om}{\Omega}
\newcommand{\om}{\omega}

\newcommand{\weak }{\, -\!\!\!\!-\!\!\!\!\rightharpoonup}
\newcommand{\weakstar }{ \overset{\, *_{\phantom{|}}}{{\smash{\weak }}\, } }

\newcommand{\larrow}{\longrightarrow}
\newcommand{\ot}{\otimes}
\newcommand{\Larrow}{\Longrightarrow}

\newcommand{\p}{\partial}
\newcommand{\sub}{\subseteq}
\newcommand{\set}{\setminus}
\newcommand{\by}{\times}

\newcommand{\tr}{\mathrm{tr}}
\DeclareMathOperator{\sgn}{sgn} 

\newcommand{\dist}{\mathrm{dist}}

\newcommand{\spn}{\mathrm{span}}

\newcommand{\bt}{\begin{theorem}}\newcommand{\et}{\end{theorem}}
\newcommand{\bd}{\begin{definition}}\newcommand{\ed}{\end{definition}}
\newcommand{\bl}{\begin{lemma}}\newcommand{\el}{\end{lemma}}
\newcommand{\beq}{\begin{equation}}\newcommand{\eeq}{\end{equation}}
\newcommand{\bc}{\begin{claim}}\newcommand{\ec}{\end{claim}}
\newcommand{\bex}{\begin{example}}\newcommand{\eex}{\end{example}}
\newcommand{\bcor}{\begin{corollary}}\newcommand{\ecor}{\end{corollary}}
\newcommand{\bp}{\begin{proof}}\newcommand{\ep}{\end{proof}}

\newcommand{\BPL}{\medskip \noindent \textbf{Proof of Lemma} }

\newcommand{\BPCOR}{\medskip \noindent \textbf{Proof of Corollary} }
\newcommand{\BPP}{\medskip \noindent \textbf{Proof of Proposition} }
\newcommand{\BPT}{\medskip \noindent \textbf{Proof of Theorem} }

\numberwithin{equation}{section}

\begin{document}

\title[Vector-valued viscosity solutions for PDE systems]{On a vector-valued generalisation of viscosity solutions for general PDE systems}
 
\author{Nikos Katzourakis}

\address{Department of Mathematics and Statistics, University of Reading, Whiteknights, PO Box 220, Reading RG6 6AX, United Kingdom}

\email{n.katzourakis@reading.ac.uk}

  \thanks{\!\!\!\!\!\!\!\!\texttt{The author has been partially financially supported by the EPSRC grant EP/N017412/1}}
  

\date{}


\begin{abstract} We propose a theory of non-differentiable solutions which applies to fully nonlinear PDE systems and extends the theory of viscosity solutions of Crandall-Ishii-Lions to the vectorial case. Our key ingredient is the discovery of a notion of extremum for maps which extends min-max and allows ``nonlinear passage of derivatives" to test maps. This new PDE approach supports certain stability and convergence results, preserving some basic features of the scalar viscosity counterpart. In this introductory work we focus on studying the analytical foundations of this new theory.
\end{abstract}

\maketitle

\tableofcontents

\section{ Introduction} \label{section0}

The theory of viscosity solutions (VS) of Crandall-Ishii-Lions is one of most successful contexts of generalised solutions in which fully nonlinear degenerate elliptic and parabolic equations can be studied effectively \cite{BCD, CESS, Ba, B, CC, CS, CGG, C2, CIL, CEL, CL, E1, E2, G, I1, IL, J1, K, Ko, L, So}. Main attributes of this approach are the flexibility in passing to limits and the strong uniqueness theorems that supports, albeit the solutions may be nowhere differentiable in any sense.

The tremendous success of VS in the last forty years in tackling a variety of problems has led many researchers to investigate possible generalisations of this theory to virtually any possible direction. Here ``generalisation" is meant either in the strict sense of either generalising the main ideas in order to apply under weaker assumptions to more equations, or in the loose sense of generalising the spirit of applicability of this approach. Without being exhaustive, some notable extensions of viscosity solution in either sense are given in \cite{BCN, Br, BM, CGG, ESp, IK, LS, LS2, S}.

It appears that the main restriction which to date has been unable to be removed is that VS apply to single equations with scalar-valued solutions, or at best to weakly coupled monotone systems which essentially can be treated component-wise as independent equations. Removing this constraint is hardly a straightforward task, as VS are essentially based on the scalar nature of the problem and on comparison arguments through the Maximum Principle and what is typically referred to as its ``calculus", by which is meant that at a point of local extremum of a twice differentiable scalar function, its gradient vanishes and its hessian has sign.

In this work we propose a vectorial generalisation of VS which applies to non-monotone fully nonlinear systems. {\it Our key ingredient is the discovery of a notion of extremum for mappings which extends min/max and allows ``nonlinear passage of derivatives" to test maps}. This notion of vectorial extremum, which is of independent interest in itself, is coined {\bf contact} and the term is inspired by the terminology of ``touching" by smooth test functions from above/below in the scalar theory of VS. {\it The notion of contact is characterised uniquely by a ``maximum principle'' type of calculus (vanishing gradient, inequality for the hessian) for vector-valued
functions}. The notion of contact of maps allows to develop a theory of generalised solutions for systems, coined {\bf contact solutions (CS)}, which supports various stability and convergence results, preserving some of the trademark attributes of VS and in some cases extending mutatis mutandis some well-known scalar facts.

Despite simple to state, the notion of contact presents unfamiliar peculiarities. The possible vectorial ``twist'' forces the notion to be \emph{functional} rather than pointwise, in the sense that extremals are \emph{maps} and not \emph{points}. Further, it is not associated with {\bf any} partial ordering. Moreover, it has \emph{order}: first order contact of two regular maps implies equality of their gradients and second order contact implies an additional certain tensor inequality for their hessians at the point of contact. Finally, it is not very easy to motivate how it arises, as its introduction is justified merely by the generalised maximum principle calculus it carries. Hence, for pedagogical reasons we have chosen to found our exposition on the conceptually equivalent notion of {\bf contact jets}, rather than contact maps. The former are sets of pointwise generalised derivatives which extend the sub/super jets of VS and vectorial extrema appear rather later in our exposition. We hope that this significantly simplifies the presentation, as contact jets bear a strong formal resemblance with their scalar counterparts in their philosophy, once one replaces the usual inequalities in the Taylor expansions by appropriate tensor inequalities. 

One of the first natural question that comes to mind of the experts of VS every time an extension of the notion appears is ``what kind of uniqueness results can be obtained". In the general vectorial case this may not be the right question to ask, but rather ``what kind of existence results can be obtained". For instance, there are striking examples of systems whose scalar counterparts have a good (existence and) uniqueness theory, whilst in the full vectorial case even determining conditions for existence is a highly non-trivial issue (see e.g.\ \cite{LO}). Hence, there are serious limitations due to vectorial obstructions and in this work there is a natural shift towards tools for existence rather than uniqueness.

We now give a brief outline of the contents of this work. In Section \ref{section1} we introduce the preliminaries of multilinear algebra required to deal with the vectorial case. In Section \ref{section2} we introduce the notion of contact jets as
sets of generalised pointwise derivatives and the notion of CS. In Section \ref{section3} we introduce the appropriate notion of
degenerate ellipticity for fully nonlinear second order systems of PDEs
that guarantee the compatibility of CS with classical solutions. In Section \ref{section4} we study the basic properties of
contact jets, relate them to classical derivatives and prove the
compatibility of CS with the classical notions. We also derive some equivalent formulations. In Section \ref{section6} we study the local structure around a point at which a contact jet exists. The main result here is a formulation of contact jets which
involves a single scalar (not tensorial) fundamental inequality. We
then deduce that obstructions arising in the vectorial case imply that
the $J^1$-theory requires a priori little H\"older space $c^{\frac{1}{2}}$
regularity and the $J^2$-theory requires a priori Lipschitz
$C^{0,1}$ regularity. In the scalar case all obstructions disappear
and the reduced theory of semi-jets $J^{2,\pm}$ applies successfully
to merely $C^0$ functions. Roughly speaking, in the general vectorial case only
``$1/2$'' of the derivatives can be interpreted weakly, the rest
``$1/2$'' must exist classically. In Section \ref{section7} we present the extremality notion of contact, establish its characterisation through a certain Calculus. Finally, in Section \ref{section8} we prove a stability result for contact solutions as an existence tool. 

We close this introduction by noting that due to considerations of size and length, in this paper we only introduce and study the abstract rudiments of CS, leaving concrete applications to specific systems for future developments.

\section{ Preliminaries, multilinear maps, tensor inequalities}
\label{section1}

\noi {\bf Preliminaries.} Let $n,N \in \N$. In what follows,
$\R^n$ will always denote the space of column vectors $\R^{n\times
1}$ equipped with the inner product $a^\top b$ and
$\R^{N\times n}$ will stand for the space of $N\times n$ matrices
equipped with inner product $P:Q:=\tr(P^\top Q)$. The respective
norms will always be the Euclidean ones $|a|=\sqrt{a^\top a}$,
$|P|=\sqrt{P:P}$. The unit open ball of $\R^n$ centred at $x$ will
be denoted by $\mB^n(x)$ and the respective unit sphere by
$\mS^{n-1}(x)$. If $x=0$, then  $\mB^n(0)$, $\mS^{n-1}(0)$ will be
denoted by $\mB^n$, $\mS^{n-1}$. If $a \in \R^n$, the sign of $a$ is
$\sgn(a):=a/|a|$ if $a\neq 0$ and $\sgn(0):=0$. The summation
convention will always be employed when repeated indices appear in a
product and repeated free indices will be denoted with a hat
$\hat{(.)}$. We henceforth reserve the letters $n,N$ for the dimensions
of $\R^n$ and $\R^N$; Greek indices $\alpha$, $\beta$, $\gamma$,...\
will always run from $1$ to $N$ and Latin $i$, $j$, $k$,...\ from
$1$ to $n$. We denote the space of linear maps $\R^n \larrow
\R^N$ by $\R^N \! \ot\R^n$ or by $\R^{N \by n}$ and  the space of linear maps $\R^{N\times n} \larrow \R^{N \times n}$ by $\R^{Nn \by Nn}$. The spaces of linear symmetric/positive endomorphisms of $\R^n,\R^N,\R^{N\by n}$ will be symbolised by the respective subscripts ``$s,+$". Let now $u : \R^n \supseteq \Om \larrow \R^N$ be a twice differentiable map. $\Om$ will always denote an open subset and ``once" or ``twice" differentiability is understood as existence of first or second order Taylor expansions. We view the gradient matrix $\D u $ and the hessian tensor $\D^2 u$ as maps
 \begin{align} \label{eq1.2a}
 \D u  \ &= \ (\D_i u_\al)e_\al \ot e_i \ \ \ \ : \ \ \R^n \supseteq \Om \larrow \R^N \!\ot
 \R^n, \\
  \D^2 u  \ &= \ (\D^2_{ij} u_\al)e_\al \ot e_{ij} \ \ : \ \ \R^n \supseteq \Om \larrow \R^N \!\ot
 \R^{n\by n}_s,
 \end{align}
 where $e_{ij}:=e_i \ot e_j$, $e_{\al i}:= e_\al \ot e_i$ and $e_\al$, $e_i$ denote the standard bases of $\R^N$ and $\R^n$ respectively. We also introduce the following \emph{contraction operation} for tensors which extends the Euclidean inner product of $\R^N \! \ot\R^n$. Let ``$\ot^{(r)}$'' denote the $r$-fold tensor product. If $S\in \ot^{(q)}\R^N \ot^{(s)} \R^n$,  $T \in \ot^{(p)}\R^N \ot^{(s)} \R^n$ and $q\geq p$, we define a tensor $S:T$ in $\ot^{(q-p)} \R^N$ by 
\beq \label{1.A}
S:T \ :=\ \big(S_{\al_q ...\al_p... \al_1 \, i_s ... i_1}  T_{\al_{p}  ... \al_1 \, i_s ... i_1}  \big) \, e_{\al_q} \ot ... \ot e_{\al_{p+1}}. 
\eeq
For example, for $s=q=2$ and $p=1$, the tensor $S:T$ of \eqref{1.A} is a vector with components $S_{\al \be i j}T_{\be ij}$ with free index $\al$ and  the indices $\be,i,j$ are contracted. In particular, in view of \eqref{1.A}, the second order linear system 
\beq \label{1.B}
A_{\al i \be j}\D^2_{ij}u_\be \, +\, B_{\al \ga k} \D_ku_\ga + C_{\al \de} u_\de\, =\, f_\al ,
\eeq
can be compactly written as $A:\D^2 u+ B:\D u +Cu=f$, where the meaning of ``$:$" in the respective dimensions is made clear by the context. Let now $\P : \R^n \larrow \R^N$ be linear map. We will \emph{always} identify linear subspaces with orthogonal projections on them. Hence, we have the split $\R^N=[\P]^\top \oplus [\P]^\bot$ where $[\P]^\top$ and $[\P]^\bot$ denote range of $P$ and nullspace of $\P^\top$ respectively. In particular, if $\xi \in \mS^{N-1}$, then $[\xi]^\top = \xi \ot \xi$ is (the projection on) the line $\spn[\xi]$ and $[\xi]^\bot$ is (the projection on) the normal hyperplane $I-\xi \ot \xi$.

\ms
\noi
{\bf Symmetrised tensor products.} The next notion plays a
crucial role to what follows. The \emph{symmetrised tensor product}
is the operation
 \beq \label{eq1.3}
\vee \ : \ \R^N \! \ot\R^N \larrow \R^{N\by N}_s \ : \ a \vee b \ := \
\frac{1}{2}\big(a \ot b + b \ot a \big).
 \eeq
Obviously, $a \vee b = b \vee a$ and $a \vee a = a \ot a$. Let us
also record the identities
\begin{align}
(a \vee b ): &\, X  \, = \, X  :(a \ot b )\, =\, a^\top X b\ ,\ \ \ X \in \R^{N\by N}_s,\label{eq1.3a}\\
|a \ot b|^2 &\, =\, |a|^2 |b|^2\ , \ \ \ |a \vee b|^2  \, =\,  \frac{1}{2}\left(|a|^2 |b|^2 + (a^\top
b)^2\right). \label{eq1.3c}
\end{align}
We will also need to consider tensor products ``$\vee$'' of higher
order between $\R^N$ and the spaces $\R^N \! \ot\R^n$ and $\R^N \!\ot
\R^{n\by n}_s$. If $\xi \in \R^N$, $\P \in \R^N \! \ot\R^n$, we view the
tensor products $\xi \ot \P$ and $\P \ot \xi$ as maps $\R^n \larrow \R^N \! \ot\R^N$. This allows to define
 \beq \label{eq1.5}
\xi \vee \P \ : \ \R^n \larrow \R^{N\by N}_s, \ \ \ \xi \vee \P \, :=\,
\frac{1}{2}\big(\xi \ot \P + \P \ot \xi \big).
 \eeq
Obviously, $ (\xi \vee \P)w= \xi \vee (\P w) $ and $\xi \vee \P = \P \vee \xi$. Similarly, if $\X=\X_{\al i
j}e_\al \ot e_{ij} \in \R^N \! \ot\R^{n\by n}_s$, we view $ \xi \ot \X $ and $\X \ot \xi$ as maps $\R^n \by \R^n \larrow \R^N \! \ot\R^N$
and we set
 \beq \label{eq1.8}
\xi \vee \X \ : \ \ \R^n \by \R^n \larrow \R^{N\by N}_s, \ \ \ \xi \vee \X\, :=\,
\frac{1}{2}\big(\xi \ot \X + \X \ot \xi \big).
 \eeq
Once again we note that $ (\xi \vee \X)(w,v)= \xi \vee (\X :w \ot v) $ and that $\xi \vee \X = \X \vee \xi$. Moreover, since
$\X_{\al i j} = \X_{\al j i}$, the tensor $\xi \vee \X$ is in
$\R^{Nn \by Nn}_s$: indeed,
 \begin{align} \label{eq1.9}
(\xi \vee \X)_{\al i\be j}\, = \, \frac{1}{2}\big(\xi_\al \X_{\be i j} + \xi_\be X_{\al i j} \big)  \ =\  (\xi \vee \X)_{\be j \al i}. 
\end{align}

\ms
\noi
{\bf Tensor Inequalities and orderings.} Let $\Xi \in \R^{Nn \by Nn}_s$. The latter space comes equipped with its natural ordering
 \beq \label{eq1.16}
\Xi \, \geq \, 0 \ \ \ \Longleftrightarrow \ \ \ \Xi : \P \ot \P \
\geq \ 0,\  \  \P \in \R^{N\by n}.
 \eeq
We now introduce a weaker notion of partial ordering in $\R^{Nn \by Nn}_s$ which emerges in the PDE theory that follows.
 \bd[Rank-One Positivity] \label{de1}
Let $\Xi \in \R^{Nn \by Nn}_s$. We say that the 4th-order tensor $\Xi$ is \emph{rank-one positive} when the quadratic form $\P \mapsto \Xi :
\P \ot \P$ is rank-one convex on $\R^N \! \ot\R^n$, that is when
 \beq \label{eq1.17}
 \eta \in \R^N, \ w \in \R^n \ \ \Larrow \ \ \Xi :(\eta \ot w) \ot (\eta \ot
 w)\, \geq \, 0.
 \eeq
In this case, we write 
\beq
\Xi \, \geq_{\ot}  0.
\eeq 
\ed
As usually, we define $\Xi \leq_{\ot} 0$ $\Leftrightarrow$ $-\Xi
\geq_{\ot} 0$ and $\Xi \leq_{\ot} \Theta$ $\Leftrightarrow$ $\Xi
- \Theta \leq_{\ot} 0$. We recall the well-known fact that the
quadratic form $\P \mapsto \Xi : \P \ot \P$ is rank-one convex
when the next function is convex on $\R$ for all $\P$, $\eta$, $w$:
 \beq \label{eq1.17a}
t \ \mapsto \ \Xi : \big(\P +t(\eta \ot w)\big) \ot \big(\P +t(\eta
\ot w)\big).
 \eeq 

\ms

We now establish that rank-one positivity ``$\geq_{\ot}$'' defines a partial ordering
in the subspace of $\R^{Nn \by Nn}_s$ consisting of
\emph{separately symmetric fourth-order tensors}:

 \bl \label{le2} Let $\Xi$, $H$, $\Theta$  be in $\R^{Nn \by Nn}_s$. Then:
\ms

 \noi (i) We have $\ \Xi \leq_{\ot} \Xi$. Also, if $\Xi \leq_{\ot} \Theta$, $\Theta \leq_{\ot} H$  then $\Xi \leq_{\ot}
 H$.
\ms

 \noi (ii) If $0 \leq_{\ot} \Xi \leq_{\ot} 0$, then
 \beq \label{eq1.18}
\big(\Xi _{\al i \be j} \, + \, \Xi _{\be j \al i}\big)e_{\al i} \ot
e_{\be j} \, = \, 0.
 \eeq
\el
\begin{corollary}[$\geq_{\ot}$ partially orderings]
The inequality of rank-one positivity defines a partial ordering in the next space of separately symmetric tensors
  \beq \label{eq1.19}
\R^{Nn \by Nn}_{s*}\  := \ \Big\{{\Xi}=\Xi_{\al i \be j}e_{\al i} \ot e_{\be
j}\ \Big| \ \Xi_{\al i \be j}= \Xi_{\be j \al i} = \Xi_{\be i \al j}
\Big\}.
 \eeq
 \end{corollary}

\BPL \ref{le2}. (i) is trivial. To see (ii), let
$\xi \in \R^N$, $w\in\R^n$. By assumption, we have $0 \leq \Xi:(\xi
\ot w)\ot(\xi \ot w) \leq 0$. Hence,
\begin{align} \label{eq1.20}
0\, = \, \Xi:(\xi \ot w)\ot(\xi \ot w)\, = \, \Xi_{\al i\be j} \xi_\al w_i \xi_\be w_j \, = \,  (\xi_\al \Xi_{\al i\be j} \xi_\be) w_i  w_j . 
\end{align}
If we set $\xi^\top \Xi \, \xi := (\xi_\al \Xi_{\al i\be j} \xi_\be)
e_i \ot e_j$, then \eqref{eq1.20} says
 \beq \label{eq1.21}
(\xi^\top \Xi \, \xi):w \ot w \ =  \ 0
 \eeq
for all $w \in \R^n$, and moreover, by the symmetries of $\Xi$, we have $\xi^\top \Xi \,
\xi \in \R^{n\by n}_s$:
\begin{align} \label{eq1.22}
(\xi^\top \Xi \, \xi)_{ij}\ &= \ \xi_\al \Xi_{\al i\be j} \xi_\be\, = \, \xi_\al \Xi_{\be j \al i} \xi_\be \
   = \ (\xi^\top \Xi \, \xi)_{ji} .
\end{align}
Hence, \eqref{eq1.21} implies for all $i,j \in \{1,...,n\}$ and
all $\xi \in \R^N$ that
 \beq \label{eq1.23}
\xi_\al \Xi_{\al i\be j}\, \xi_\be = 0.
 \eeq
By interchanging in \eqref{eq1.23} $i$ and $j$ and employing
that $\Xi_{\al i\be j}=\Xi_{\be j \al i}$, we have
 \beq \label{eq1.24}
\xi_\al \Xi_{\be i \al j}\, \xi_\be \, = \, 0.
 \eeq
By \eqref{eq1.23} and \eqref{eq1.24}, for all $i,j$ fixed we have
 \beq \label{eq1.25}
\xi_\al \big(\Xi_{\al i\be j}\, + \, \Xi_{\be i\al
j}\big) \xi_\be \, = \, 0.
 \eeq
Since $\big(\Xi_{\al i\be j} +  \Xi_{\be i\al j}\big) e_\al \ot e_\be$ belongs to $\R^{N\by N}_s$, we obtain \eqref{eq1.18} as
desired. \qed

\ms

\begin{remark} 
It is evident that $\R^{Nn \by Nn}_{s*}$ is a proper subspace of $\R^{Nn \by Nn}_s$
and that it can be equipped with both partial orderings ``$\leq$''
and ``$\leq_{\ot}$''. It also evident that ``$\leq$'' is a stronger
notion than ``$\leq_{\ot}$'', in the sense that $ \Xi \geq 0$ implies $ \Xi \geq_{\ot}
 0$. The known examples of rank-one convex quadratic form which are not convex imply that rank-one positivity is genuinely weaker that positivity.
\end{remark}

\section{ Contact solutions for fully nonlinear PDE systems}
\label{section2}

In this section we introduce the basics of a theory of
non-differentiable solutions which applies to fully nonlinear
systems of partial differential equations of the form
 \beq \label{3.1}
\F(\cdot , u,\D u ,\D^2 u   )\, = \, 0 ,
  \eeq
where $u:\R^n \supseteq \Om \larrow \R^N$ and
 \beq \label{3.2}
\F \  :\  \ \Om \by \R^N \by (\R^N \! \ot\R^n) \by (\R^N \! \ot\R^{n\by n}_s) \larrow \R^N.
  \eeq
The arguments of the nonlinearity $\F$ will be denoted by $\F\big(x,\eta,\P,\X
\big)$. For the moment, the only assumption that needs to be imposed to $\F$ is
mere local boundedness. Hence, we allow for \emph{discontinuous coefficients}. Later we will assume continuity and an appropriate notion of ellipticity, in order to assure compatibility of generalised and classical solutions. Our
notion of solution allows to interpret merely continuous maps as solutions to the PDE system \eqref{3.1}. The point of view is to relax $\D u $ and $\D^2 u $ to certain
generalised \emph{pointwise} derivatives and relax equality in
\eqref{3.1} to appropriate inequalities, when $\F$ is evaluated at these generalised derivatives.

\bd[Contact jets] \label{de2.1} Let $u : \R^n \supseteq \overline{\Om} \larrow \R^N$ be a continuous map, $x\in \overline{\Om}$  and $\xi \in \mS^{N-1}$. The \emph{first contact $\xi$-jet of $u$ at $x$} is the set
of generalised derivatives
 \begin{align} \label{eq2.3}
J^{1,\xi}u(x)  \  := \  \Big\{\ & \P \in \R^N \! \ot\R^n \ \big| \text{ as \ } \overline{\Om} \ni z\rightarrow x,\nonumber\\
 &\xi \vee \big[u(z)-u(x)-\P (z-x)\big] \leq o(|z-x|) \ \Big\}.
 \end{align}
The \emph{second contact $\xi$-jet of $u$ at $x$} is the set of
generalised derivatives
 \begin{align} \label{eq2.4}
J^{2,\xi}u(x)  \  :=& \  \Bigg\{ \ (\P,\X) \in \R^N \! \ot (\R^N \by \R^{n\by n}_s)\ \Big|\text{ as \ } \overline{\Om} \ni z\rightarrow x, \nonumber\\
& \xi \vee \Big[u(z)-u(x)-\P (z-x)  -\frac{1}{2}\X :(z-x) \ot (z-x) \Big] \leq o(|z-x|^2) \ \Bigg\}. 
 \end{align}
 \ed

\begin{remark}
The meaning of ``$o(1)$'' in \eqref{eq2.3}, \eqref{eq2.4} is that
there exists a continuous matrix-valued map $T : \R^n \set \{0\}
\larrow \R^{N\by N}_s$ such that $|T(y)| \rightarrow 0$ as $y\rightarrow 0$. The meaning of the rest quantities appearing is that given in formulas \eqref{eq1.3}-\eqref{eq1.8}. In particular, matrix inequalities are considered in $\R^{N\by N}_s$. The necessity to define generalised derivatives on boundary points $x\in \p \Om$ of closed sets $\overline{\Om}$ stems from the necessity to consider boundary value problems for PDE systems, but also for technical reasons arising in our subsequent analysis. If $x \in \Om$, since $\Om$ is open the statement ``$\overline{\Om} \ni z$" which means ``convergence in $\overline{\Om}$" can be dropped. 

In the scalar case of $N=1$, we have $\mS^0=\{-1,+1\}$ and
$J^{1,\xi}$,$J^{2,\xi}$ reduce to the semi-jets $J^{1,\pm}$ and $J^{2,\pm}$  of  VS. Indeed, $(\pm 1)\vee a = (\pm 1) \ot a = \pm a$ for any $a \in \R$ and \eqref{eq2.4} reduces to inequality in $\R$. Moreover, by applying ``$:(\xi \ot \xi)$'' to \eqref{eq2.4} we deduce that the ``scalar'' $\xi$-projection $\xi^\top u$ of $u$ along $\xi \ot \xi$ satisfies  $(\xi^\top \P, \xi^\top \X)\in J^{2,\pm}(\xi^\top u)(x)$, and similarly for \eqref{eq2.3} (see \eqref{eq1.3a}). 
\end{remark}

In the following we will also need to consider closures of contact jets:

\bd[Contact jet closures] \label{de2.1b} Let $u : \R^n \supseteq \overline{\Om} \larrow \R^N$ be a continuous map, $x\in \overline{\Om}$  and $\xi \in \mS^{N-1}$. The \emph{first contact $\xi$-jet closure of $u$ at $x$} is 
 \begin{align} \label{eq2.3a}
\overline{J}^{1,\xi}u(x)  \  := \  \Big\{\ & \P \in \R^N \! \ot\R^n \ \big| \ \exists\  (\xi_m,x_m,\P_m,)\to (\xi,x,\P) \nonumber\\
& \text{ as }m\to \infty \ : \ \P_m \in J^{1,\xi_m}u(x_m)\Big\}.
 \end{align}
The \emph{second contact $\xi$-jet closure of $u$ at $x$} is
 \begin{align} \label{eq2.4a}
\overline{J}^{2,\xi}u(x)  \  := \  \Big\{\ & \P \in \R^N \! \ot(\R^n \by \R^{n\by n}_s) \ \big| \ \exists\ (\xi_m,x_m,\P_m,\X_m)\to \nonumber\\
& (\xi,x,\P,\X)  \text{ as }m\to \infty \ : \ (\P_m,\X_m) \in J^{2,\xi_m}u(x_m)\Big\}.
 \end{align}
 \ed
The main difference of \eqref{eq2.3a}, \eqref{eq2.4a} compared to their scalar counterparts is that we approximate in the direction $\xi$ as well. When $N=1$ no such option is available, since $\mS^0=\{-1,+1\}$ is totally disconnected. Before giving our notion of solution, we need one more definition, which we state only for the second order case. 

\bd[Envelopes of discontinuous coefficients] Given $\xi \in \mS^{N-1}$ and consider the map $\F$ of \eqref{3.2} which we assume is locally bounded. The \emph{$\xi$-Envelope $\xi^*F$ of $\F$} is the upper semi-continuous envelope of the projection $\xi^\top \F$:
\begin{align}
\xi^*\F(x,\eta, \P,\X)\ :=&\ \lim_{\e \to 0}\sup\Big\{\xi^\top \F(y,\theta, \Q,\Y) \ :\ |x-y|\\
& \ \ \  +|\eta-\theta|+|\P-\Q|+|\X-\Y|\leq \e\Big\}. \nonumber
\end{align}
\ed

We now proceed to the main notions of solutions we will use in this work.

\bd[Contact solutions for second order systems] \label{de2.3}

\noi Consider the map $\F$ of \eqref{3.2} and suppose it is locally bounded. The continuous map $u :\R^n \supseteq \Om \larrow \R^N$ is called a \emph{contact solution to} \eqref{3.1} on $\Om$ when for any $x\in \Om$ and $\xi \in \mS^{N-1}$ we have
 \begin{align} \label{eq2.8}
(\P,\X) \in \overline{J}^{2,\xi}u(x) \ \ \Larrow \ \ \xi^*
\F\big(x,u(x),\P,\X\big) \, \geq \, 0.
 \end{align}
 \ed
Similarly, one can specialise the notion for first order systems as follows.

\bd[Contact solutions for first order systems] \label{de2.2}

\noi Suppose the map
\beq
\F \ : \ \Om \by \R^N \by (\R^N \! \ot\R^n) \larrow \R^N
\eeq
is locally bounded. The continuous map $u :\R^n \supseteq \Om \larrow \R^N$ is called a
\emph{contact solution to}
 \beq \label{eq2.5}
\F(\cdot , u,\D u )\, = \, 0
  \eeq
on $\Om$, when for all $x\in \Om$ and $\xi \in \mS^{N-1}$ we have
 \begin{align} \label{eq2.6}
\P \in \overline{J}^{1,\xi}u(x) \ \ \Larrow \ \ \xi^* \F\big(x,u(x),\P\big) \,
\geq \, 0.
 \end{align}
 \ed

\begin{remark} We observe that for $N=1$, CS reduce  to viscosity
solutions (up to a difference in the sign convention in the
inequality). However, the new ingredients in the contact notions which are \emph{not} component-wise will lead to genuinely vectorial phenomena. 
\end{remark}

The new objects $J^{1,\xi}$, $J^{2,\xi}$ will be studied thoroughly
later. Before presenting some explicit calculations of contact jets for a typical map to illustrate the working philosophy (which is analogous to the scalar case), we present a reformulation of Definitions \ref{de2.3}-\ref{de2.2}.

\bl[Alternative definitions]  \label{l1a} In the setting of Definitions \ref{de2.3}-\ref{de2.2}, the implications \eqref{eq2.6}, \eqref{eq2.8} can be respectively replaced by
 \begin{align} \label{eq2.6a}
\P \in J^{1,\xi}u(x) \ \ &\Larrow \ \ \xi^*
\F\big(x,u(x),\P\big) \, \geq \, 0,\\
\label{eq2.8a}
(\P,\X) \in J^{2,\xi}u(x) \ \ &\Larrow \ \ \xi^* \F\big(x,u(x),\P,\X\big) \, \geq \, 0.
 \end{align}
If moreover the nonlinearity $\F$ is continuous, we can replace the $\xi^*F$ by $\xi^\top \F$.
\el

\BPL \ref{l1a}. For brevity we exhibit only the second order case. Obviously, $ J^{2,\xi}u(x)  \sub  \overline{J}^{2,\xi}u(x)$. Conversely, assume \eqref{eq2.8a} and fix $(\P,\X) \in \overline{J}^{2,\xi}u(x)$. Then, there is a sequence $(\xi_m,x_m,\P_m,\X_m)\to  (\xi,x,\P,\X)$ as $m\to \infty$ and also $ (\P_m,\X_m) \in J^{2,\xi_m}u(x_m)$. By \eqref{eq2.8a}  and since $u$ is continuous and $\F$ is locally bounded, there exists a bounded open set $B$ centred at $(\xi,x,\P,\X)$ such that
\begin{align}
0\, & \leq \, (\xi_m)^*\F\big(x_m,u(x_m),\P_m,\X_m\big)\\
     &\leq \, \xi^*\F\big(x_m,u(x_m),\P_m,\X_m\big)\ +\ |\xi-\xi_m|\big(\sup_B |\F|\big) \nonumber
\end{align}
for $m$ large. Since $\xi^*F$ is upper semi-continuous, by letting  $m\to \infty$  we obtain that $\xi^*\F(x,u(x),\P,\X)\geq 0$. Hence, the map $u$ is a contact solution.  
\qed

\section{ Basic calculus of generalised derivatives}
\label{section4}

In this section we examine the pointwise generalised derivatives $ J^{1,\xi}u(x)$ and $ J^{2,\xi}u(x)$ and relate them with the classical ones $\D u (x),\D^2 u   (x)$. We begin with two simple algebraic results will turn out to be essential tools.

  \bl[Spectral decomposition of symmetrised tensor products] \label{le1}

\noi Let $R \in \R^N$ and $\xi \in \mS^{N-1}$. Then, $\xi \vee R$ is a symmetric matrix with rank at most $2$ and its the spectrum of consists of at most three distinct eigenvalues $\la^- \leq 0 \leq \la^+$, given by
 \beq \label{eq1.10}
\si (\xi \vee R) \, = \, \left\{-\frac{1}{2}\left(|R|-\xi^\top
R\right),\ 0,\ \frac{1}{2}\big(|R|+\xi^\top R\big) \right\}.
 \eeq
The respective eigenspaces are
 \beq \label{eq1.11}
\mathrm N\big(\xi \vee R \ -\ \la I\big)\ = \
 \left\{\begin{array}{l}
\ \spn[\xi - \sgn(R)], \ \ \ \la = \la^- ,\\
\left(\spn[\{\xi,R\}]\right)^\bot, \hspace{18pt} \la = 0, \\
\ \spn[\xi + \sgn(R)] , \ \ \ \la = \la^+ .
 \end{array}
 \right.
 \eeq
  \el

\BPL \ref{le1}.  If $R=0$ or $R$ is co-linear to $\xi$, the result is obvious. If $R,\xi$ are linearly independent, we observe that $\mathrm N\big(\xi \vee R\big)=
\left(\spn[\{\xi,R\}]\right)^\bot$: indeed, for all $\eta \in\R^N$,
we have the identity
 \beq \label{eq1.13}
 (\xi \vee R)\eta \, = \, \left(\frac{R^\top \eta}{2}\right)\xi \ + \
 \left(\frac{\xi^\top \eta}{2}\right)R.
 \eeq
Hence, $\eta$ is normal to both $\xi$ and $R$ if and only if $(\xi
\vee R)\eta =0$. By the Spectral Theorem, $\xi \vee R$ has at most
three distinct eigenvalues $\la^-$,$0$, $\la^+$ and
 \beq \label{eq1.14}
\mathrm N\big(\xi \vee R \ -\ \la^- I\big)  \oplus \mathrm N\big(\xi \vee R \ -\
\la^+ I\big) \, = \, \spn[\{\xi,R\}].
 \eeq
We now employ \eqref{eq1.13} to check directly that
 \beq \label{eq1.15}
(\xi \vee R)\left(\xi \pm \frac{R}{|R|}\right) \, = \, \la^\pm
\left(\xi \pm \frac{R}{|R|}\right)
 \eeq
with $\la^\pm$ as in \eqref{eq1.10}. The lemma follows.
\qed
\ms

We now show that symmetric products $\xi \vee (\cdot)$ coupled by the inequality $\geq_{\ot}$ induce ``directed'' orderings.

 \begin{proposition}[Induced partial orderings] \label{pr1} Let $\xi$ be in $\mS^{N-1}$ and $\xi^\bot = I - \xi \ot \xi$.

\ms
\noi (i) If $v \in \R^N$, then
 \begin{align} \label{eq1.30}
\xi \vee v \, \leq \, 0 \ \Leftrightarrow \ 
\left\{
\begin{array}{l}
v \, =\, (\xi^\top v)\xi \\
\xi^\top v \, \leq \, 0
\end{array}
\right.
 \ \Leftrightarrow \   
\left\{
\begin{array}{l}
\xi^\bot v \, =\, 0 \\
\xi^\top v \, \leq \, 0
\end{array}
\right.
 \  \Leftrightarrow \  v \, = \, -|v|\xi.
 \end{align}
\ms

\noi (ii) If $\X \in \R^N \! \ot\R^{n\by n}_s$, then
 \begin{align}  \label{eq1.32}
\xi \vee \X  \, \leq_{\ot}  0 \ 
\ \Leftrightarrow \ \   
\left\{
\begin{array}{l}
\X \, =\,\xi \ot  (\xi^\top \X) \\
\xi^\top \X \, \leq \, 0
\end{array}
\right.
\Leftrightarrow \  
\left\{
 \begin{array}{l}
\xi^\bot \X \, =\, 0 \\
\xi^\top \X \, \leq \, 0
\end{array}
\right. 
 \Leftrightarrow \ \xi \vee \X \, \leq \, 0  . \nonumber
 \end{align}
  \end{proposition}

\noi In particular, it follows that the orderings $\leq$ and $\leq_{\ot}$ coincide on the
cone
 \beq
 \Big\{\eta \vee \Y \
\Big| \ \eta \in \R^N, \ \Y \in \R^N \! \ot\R^{n\by n}_s \Big\}
 \eeq
which is a subspace of the space \eqref{eq1.19} of separately
symmetric tensors.

\BPP \ref{pr1}. (i) By Lemma \ref{le1}, $\xi \vee v \leq 0$ if and only if $\max \si (\xi \vee v) \leq 0$, hence if and only if $\frac{1}{2}(|v| +\xi^\top v)=0$ and this says $v = -|v|\xi$.
The latter is equivalent to $v = (\xi^\top v)\xi$ with $\xi^\top v
\leq 0$ and to $\xi^\bot v=0$ with $\xi^\top v\leq 0$.

\noi (ii)  Suppose that $\xi \vee
\X \leq_{\ot} 0$ and fix $\eta \in \R^N$ and $w \in \R^n$. Then, we
have
 \begin{align} \label{eq1.33}
0 \ & \geq \ (\xi \vee \X ):(\eta \ot w) \ot (\eta \ot w)
\nonumber\\
& =\ \frac{1}{2}\Big(\xi_\al \X_{\be i j } \, + \, \xi_\be \X_{\al i
j }\Big) \eta_\al w_i \eta_\be w_j \\
& =\ \frac{1}{2}\left[\big(\xi_\al \X_{\be i j }w_i  w_j \big) \ + \
\big( \xi_\be \X_{\al i j } w_i  w_j \big) \right] \eta_\al \eta_\be \nonumber\\
& = \ \big(\xi \vee (\X:w \ot w )\big) :\eta \ot \eta . \nonumber
 \end{align}
By \eqref{eq1.33}, we obtain for any $w$ fixed that $\xi \vee (\X :w \ot w) \leq 0$. By employing (i) to the vector $v:=\X:w \ot w$, we see that
 \[
 \begin{split}
 \xi^\top (\X : w \ot w) \leq  0\ ,\ \ \ \ 
 \big(\X - \xi \ot(\xi^\top \X) \big) : w \ot w  =  0,
 \end{split}
 \]
for any $w$ fixed. Since $w$ is arbitrary, we obtain the desired decomposition which can be recast as $\xi^\bot\X=0$, $\xi^\top\X \leq 0$. Finally, by assuming the latter decomposititon and fixing $\P \in \R^N \! \ot\R^n$, we have
 \begin{align}
(\xi \vee \X) : \P \ot \P \  & =\ \big(\xi \vee \xi
\ot (\xi^\top \X)\big) : \P \ot \P \nonumber\\
& =\ \xi_\al \xi_\be \xi_\gamma \X_{\gamma i j} \P_{\al i} \P_{\be
j} \nonumber\\
& =\ (\xi_\gamma \X_{\gamma i j})(\xi_\al \P_{\al i}) (\xi_\be
\P_{\be j})\\
& = \ (\xi^\top \X) : (\xi^\top \P) \ot (\xi^\top \P)  \nonumber\\
& \leq \ 0 \nonumber
 \end{align}
and the last inequality follows by $\xi^\top \X \leq 0$. Hence, $\xi \vee \X \leq 0$
as desired. Finally, the implication $\xi \vee \X  \leq  0 \Rightarrow  \xi \vee
\X  \leq_{\ot}  0 $ is trivial. \qed

\ms

Now we relate generalised and classical pointwise derivatives.

 \bt[Contact jets and derivatives]  \label{th2} Let $u : \R^n \supseteq \Om \larrow \R^N$ be a map which is continuous at $x\in \Om$.

\noi (a) If there exists \emph{one direction} $\xi \in \mS^{N-1}$
such that both $J^{1,\pm \xi}u(x)$ are nonempty, then $u$ is
differentiable at $x$ and both $J^{1, \pm \xi}u(x)$ are
singletons with element the gradient:
 \beq \label{eq4.3}
J^{1,\pm \xi}u(x) \neq \emptyset \ \ \ \Larrow \ \ \  J^{1, \xi}u(x) =
J^{1, -\xi}u(x) = \big\{\D u (x)\big\}.
 \eeq

\noi (b) If $u$ is differentiable at $x$, then for all $\xi \in
\mS^{N-1}$ the sets $J^{1,\xi}u(x)$ are singletons with element the gradient:
 \beq \label{eq4.4}
 J^{1,\xi}u(x)  \,=\,  \big\{\D u (x)\big\}.
 \eeq
Moreover, whenever $(\D u (x),\X^\pm) \in J^{2,\pm \xi}u(x)\neq \emptyset$, we have the
inequality 
\beq
\xi \vee \big[\X^- - \X^+ \big]\, \leq_{\ot} 0 
\eeq
which is equivalent to
 \beq \label{eq4.6}
 \xi^\bot \big[\X^- - \X^+\big]\, =\, 0\ , \ \ \ \xi^\top \big[\X^- - \X^+\big]\, \leq \, 0.
 \eeq

\noi (c) If $u$ is twice differentiable at $x$, then for all $\xi
\in \mS^{N-1}$ the sets $J^{2,\xi}u(x)$ are nonempty, they
contain $(\D u (x),\D^2 u   (x))$ and also
 \begin{align} \label{eq4.5a}
 J^{2,\xi}u(x) \, =\,  \Big\{\big(\D u (x),\D^2 u   (x) +\xi \ot A \big)\ :  \ A\geq 0\  \Big\} .
 \end{align}
Moreover, we have the characterisations 
 \begin{align} \label{eq4.5}
 J^{2,\xi}u(x)     & =  \Big\{(\D u (x),\X) \ : \ \xi \vee
  \big[\D^2 u   (x) - \X\big]\leq_\ot 0 \ \ \Big\} \nonumber\\
 & =  \Big\{(\D u (x),\X) \ : \ \xi \vee
  \big[\D^2 u   (x) - \X\big]\leq 0 \ \ \Big\} \\
  & =  \Bigg\{(\D u (x),\X) \  
: \bigg\{
\begin{array}{l}
\xi^\top \big[\D^2 u   (x) - \X\big]\leq 0,\ms\\ 
\xi^\bot  \big[\D^2 u   (x) - \X\big]= 0 
\end{array}
\Bigg\}.\nonumber
 \end{align}

\noi (d) If $v : \R^n \supseteq \Om \larrow \R^N$ is twice differentiable at $x$ and $\la,\mu \geq 0$, then
\beq
J^{2,\xi}(\la u +\mu v)(x)\, =\, \la J^{2,\xi}u(x)\, +\, \mu \big(Dv(x),D^2v(x)\big).
\eeq
 \et

\BPT \ref{th2}. (a) Let $\P^\pm \in J^{1,\pm \xi}u(x) \neq \emptyset$.
Then, by \eqref{eq2.3}, we have
 \begin{align}
 \pm\xi \vee \big[u(z+x)-u(x)-\P^\pm z\big]:\eta \ot \eta\ \leq\ o(|z|),  \label{eq4.7} 
 \end{align}
if $\eta \in \R^N$, where ``$o(1)$'' is realised by $T(z): \eta \ot \eta$. We set $z:=\e w$ for $\e>0$ and $|w|=1$ fixed and add the $\pm$ inequalities in \eqref{eq4.7} to
obtain
 \beq \label{eq4.9}
 \xi \vee \big[(\P^- -\P^+)w\big]:\eta \ot \eta \, \leq\, \frac{o(\e)}{\e}
 \eeq
as $\e \rightarrow 0^+$. By taking the limit to \eqref{eq4.9} and
then replacing $w$ with $-w$, we obtain $\xi \vee \big[(\P^+
-\P^-)w\big]:\eta \ot \eta =0$. By applying Proposition \ref{pr1},
we find that $\xi \vee (\P^- -\P^+)$ vanishes. By Lemma \ref{le1}, we obtain that zero is the
unique eigenvalue of $\xi \vee (\P^+ -\P^-)$, that is, $\si\big(\xi
\vee (\P^+ -\P^-)\big)=\{0\}$. Hence, we have
 \beq
\big|(\P^+ - \P^-)w\big| \ =\ \pm \xi^\top \big((\P^+ - \P^-)w\big)\
=\ 0,
 \eeq
for any $w \in \mS^{n-1}$. As a result we get $\P^+ = \P^-$. Let us denote
their common value by $\P$. Then, by \eqref{eq4.7} we
have
 \beq \label{eq4.10}
\xi \vee \big[u(z+x)-u(x)-\P z\big]:\eta \ot \eta\ =\ o(|z|),
  \eeq
as $z\rightarrow 0$. Since the numerical radius (cf.\
\cite{Lax})
 \beq \label{eq4.11}
\|A\|\ := \ \max_{\eta \in \mS^{N-1}}\big|A:\eta \ot \eta\big|
 \eeq
is a norm on $\R^{N\by N}_s$ equivalent to the Euclidean, \eqref{eq4.10}
implies as $z \rightarrow 0$ that
 \begin{align} \label{eq4.12}
 o(|z|)\ &= \ \max_{\eta \in \mS^{N-1}}
 \big|\xi \vee [u(z+x)-u(x)-\P z ]:\eta \ot \eta \big| 
 \\
& \geq \ \frac{1}{C}\big|\xi \vee [ u(z+x)-u(x)-\P z ] \big|,  \nonumber
 \end{align}
for some $C>0$. By applying \eqref{eq1.3c}, \eqref{eq4.12} gives
 \[
 \begin{split} \label{eq4.13}
 o(|z|)\ & = \ \frac{1}{2C^2} \Big\{ \big|u(z+x)-u(x)-\P z \big|^2 \,  + \, \big[\xi^\top ( u(z+x)-u(x)-\P z ) \big]^2 \Big\}  \\
& \geq \ \frac{1}{2C^2}\big| u(z+x)-u(x)-\P z \big|^2,  
 \end{split}
 \]
as $z \rightarrow 0$. Consequently, we have $\P=\D u (x)$ and
$J^{1,\xi}u(x)=\big\{\D u (x)\big\}$.

\medskip

\noi (b) If $u$ is differentiable at $x$, by applying $\xi \vee
(\cdot )$ to the Taylor expansion $u(z)-u(x)-\D u (x)z\, = \, o(|z|)$
which holds as $z \rightarrow 0$, we discover $\big\{\D u (x)\big\} \sub J^{1,\xi}u(x)$, for any $\xi \in \mS^{N-1}$. Since then both $J^{1,\pm \xi}u(x) \neq \emptyset$, application of (a)
implies that $J^{1,\xi}u(x)=\big\{\D u (x)\big\}$.  Let $(\P,\X) \in
J^{2,\xi}u(x)$. By \eqref{eq2.4}, we have
 \begin{align} \label{eq4.14}
\xi \vee \big[u(z+x)-u(x)-\P z \big] \, \leq\, \frac{1}{2}\, \xi \vee \big[\X:z \ot z\big] \, + \, o(|z|^2)\,= \, o(|z|), 
 \end{align}
as $z \rightarrow 0$. Hence, (a) implies that $\P=\D u (x)$ whenever
$(\P,\X) \in J^{2,\xi}u(x)$ and $u$ is differentiable. If  $\X^\pm \in J^{2,\pm \xi}u(x)$, we have as $z\rightarrow 0$
 \begin{align}
 & \pm \xi \vee\Big[u(z+x)-u(x)-\D u (x)z -\frac{1}{2}\X^\pm :z\ot
z\Big] \ \leq\ o(|z|^2),  \label{eq4.15}
 \end{align}
We set $z:=\e w$ for $\e>0$, $|w|=1$ and add the $\pm$ inequalities in \eqref{eq4.15} to
find
 \beq \label{eq4.17}
\xi \vee\Big[\big(\X^-  -  \X^+ \big):w\ot w \Big]:\eta \ot \eta \
\leq \ \frac{o(\e^2)}{\e^2},
 \eeq
as $\e\rightarrow 0^+$, for all $\eta \in \R^N$. By passing to the
limit in \eqref{eq4.17} we obtain
 \begin{align} \label{eq4.18}
 0 \ & \geq \ \xi \vee\Big[\big(\X^-  -  \X^+ \big):w\ot w \Big]:\eta \ot \eta  \nonumber\\
     & = \ \frac{1}{2}\xi_\al \big(\X^-  -  \X^+ \big)_{\be i j} w_i
    w_j \eta_\al \eta_\be \ +
    \ \frac{1}{2}\xi_\be \big(\X^-  -  \X^+ \big)_{\al i j} w_i w_j \eta_\al \eta_\be \nonumber\\
    & = \ \frac{1}{2} \left[\xi_\al \big(\X^-  -  \X^+ \big)_{\be i j} \ +
    \ \xi_\be \big(\X^-  -  \X^+ \big)_{\al i j} \right](\eta_\al w_i)(\eta_\be w_j) \nonumber\\
     & = \ \Big[\xi \vee \big(\X^-  -  \X^+ \big)\Big]:(\eta \ot w) \ot
     (\eta \ot w), \nonumber
     \end{align}
for all $\eta \in \R^N$, $w \in \R^n$. Hence, by Definition
\ref{de1} we obtain $\xi \vee \big(\X^-  -  \X^+ \big)\leq_\ot 0$. By Proposition \ref{pr1}, the equivalence of the rank-one inequality with \eqref{eq4.6} follows.

\medskip

\noi (c) We first observe that by applying Proposition \ref{pr1}, all four sets appearing in the right hand sides of of \eqref{eq4.5a}, \eqref{eq4.5}  are equal. Hence, it suffices to prove that $J^{2,\xi}u(x)$ equals one of those. By applying $\xi \vee (\cdot )$ to the Taylor expansion
 \beq \label{eq4.14}
u(z+x)-u(x)-\D u (x)z-\frac{1}{2}\D^2 u   (x):z\ot z\,= \, o(|z|^2),
 \eeq
which holds as $z \rightarrow 0$, we find that $(\D u (x),\D^2 u   (x)) \in
J^{2,\xi}u(x)$ for all $\xi \in \mS^{N-1}$. By applying (b) for $\X^- := \D^2 u   (x)$ and $\X^+ := \X$, we obtain the inclusion
  \beq
\label{eq4.19}
 J^{2,\xi}u(x) \ \sub \ \Big\{(\D u (x),\X)\ \Big| \ \xi \vee
  \big[\D^2 u   (x) - \X\big]\leq_{\ot} 0\ \Big\}.
 \eeq
For the reverse inclusion, let us assume that $\xi \vee \big[\D^2 u   (x)
- \X\big]\leq_{\ot} 0$. Then by applying $\xi \vee (\cdot )$ to \eqref{eq4.14}, we have
 \begin{align} \label{eq4.20}
\xi \vee\Big[u(z+x)- & u(x)-  \D u (x)z\Big]  :\eta \ot \eta \nonumber\\
& = \, \frac{1}{2}\big[\xi \vee \D^2 u   (x)\big]:\big(\eta \ot z\big) \ot \big(\eta \ot z\big) \, + \, o(|z|^2) \\
& \leq \, \frac{1}{2}\big[\xi \vee \X\big]:\big(\eta \ot z\big)
\ot \big(\eta \ot z\big)\, + \, o(|z|^2),  \nonumber
\end{align}
as $z \rightarrow 0$. By \eqref{eq4.20}, we obtain
 \beq \label{eq4.21}
\xi \vee\Big[u(z+x)-u(x)- \D u (x)z - \frac{1}{2}\X:z\ot z\Big] \, \leq\, o(|z|^2) ,
 \eeq
as $z \rightarrow 0$. Hence, $(\D u (x),\X)\in
J^{2,\xi}u(x)$.

\medskip

\noi (d) follows easily by arguing similarly as in (a), (b), (c).     \qed

\medskip

The next lemma is an essentially scalar fact which will allow to
formulate equivalent definitions of $J^1$, $J^2$. For the proof we refer to \cite{K}.

\bl \label{le4.3} Suppose $T : \R^n  \larrow \R^{N\by N}_s$ is a continuous symmetric
tensor map satisfying $|T(z)| \rightarrow 0$ as $z\to 0$.
Then, there exists an increasing function $\tau \in C^2(0,+\infty)$
with $\tau(0^+)=0$ such that $T(z) \leq \tau(|z|) I$, as $z\to 0$.
 \el

Now we derive equivalent formulations of contact jets. We shall consider only the case on $J^2$; analogous results hold for  $J^1$, with the obvious modifications. For simplicity we fix $x=0$.

\begin{proposition}[Equivalent formulations of $J^{2}$] \label{th3} Suppose $u
: \R^n \supseteq \overline{\Om} \larrow \R^N$ is continuous at $0\in \overline{\Om}$ and let $\xi \in \mS^{N-1}$ and $(\P,\X)\in \R^N \! \ot(\R^n \by \R^{n\by n}_s)$.  The following are equivalent:

\ms

 \noi (a) $(\P,\X) \in J^{2,\xi}u(0)$.

\noi (b) There exists $ \tau \in C^2(0,\infty)$ increasing with
$\tau(0^+)=0$ such that, as $\overline{\Om} \ni z\rightarrow 0$
 \beq \label{eq4.28}
 \xi \vee \left[u(z)- u(0) -\P z - \frac{1}{2}\X :z \ot z \right]\leq \ \tau(|z|)|z|^2I,
 \eeq

\noi (c) We have 
 \beq \label{eq4.29}
\underset{\overline{\Om} \ni z\to 0}{\overline{\lim}}\max_{|\eta|=1} \frac{1}{|z|^2}\left( \xi \vee \left[u(z)- u(0)- \P z -  \frac{1}{2}\X :z \ot z\right]\right):\eta \ot \eta \ \leq\ 0.
 \eeq

\noi (d)  We have as $\overline{\Om} \ni z\to 0$ that
 \beq \label{eq4.30}
 \max \si \left(  \xi \vee\left[u(z)- u(0) -\P z -  \frac{1}{2}\X :z \ot z \right] \right)\, = \, o(|z|^2).
 \eeq

\noi (e)  We have as $\overline{\Om} \ni z \to 0$ that
 \begin{align} \label{eq4.31}
&\left|u(z)- u(0)- \P z -  \frac{1}{2}\X :z \ot z \right| \nonumber\\
 &\ \  \ \ +\xi^\top \! \left(u(z)- u(0)- \P z -  \frac{1}{2}\X :z \ot z\right)\, = \, o(|z|^2).
 \end{align}
 \end{proposition}

Note that for $N=1$ we recover known properties of scalar semi-jets $J^{2,\pm}$. In particular, $(e)$ reduces to differentiability of the positive part $\big(u(z)- u(0)- \P z -  \frac{1}{2}\X$ $:z \ot z\big)^+$. The proof of proposition \ref{th3} is very simple and therefore we omit it.

\medskip

The following simple properties of $J^{1,\xi}$, $J^{2,\xi}$ are in
complete analogy with the scalar counterparts $J^{1,\pm}$,
$J^{2,\pm}$ and are a direct consequence of Proposition \ref{th3}.

\begin{proposition} \label{pr2} Let $u : \R^n \supseteq \Om \larrow \R^N$ be continuous at $x\in \Om$ and let $\xi \in \mS^{N-1}$.

\noi (a) Both $J^{1,\xi}u(x)$, $J^{2,\xi}u(x)$ are convex subsets of
$\R^N \! \ot\R^n$ and $\R^N \! \ot\R^{n\by n}_s$ respectively.

\noi (b) $J^{1,\xi}u(x)$ is closed in $\R^N \! \ot\R^n$. Moreover, for
any $\P \in \R^N \! \ot\R^n$, the ``slice"
 \beq \label{eq4.40}
\Big\{ \X \in \R^N \! \ot\R^{n\by n}_s \ \big| \ (\P,\X) \in J^{2,\xi}u(x)
\Big\}
 \eeq
is closed in $\R^N \! \ot\R^{n\by n}_s$.

\noi (c) If $J^{2,\xi}u(x) \neq \emptyset$, then it has infinite diameter. Moreover,
 \beq \label{eq4.41}
(\P,\X) \in J^{2,\xi}u(x) \ \ \ \Larrow \ \ \ (\P,\X+\xi \ot A) \in
J^{2,\xi}u(x)\ , \ \ A\geq 0.
 \eeq
\end{proposition}

\BPP \ref{pr2}. (a) is obvious. For (b), it suffices to establish
that the set \eqref{eq4.40} is closed, since the other is similar.
Let $(P,\X_m) \in J^{2,\xi}u(x)$ and $\X_m \larrow \X_\infty$ as
$m\rightarrow \infty$. Fix $\e >0$. Then, there is an $m(\e)\in \N$
such that $|\X_{m(\e)} -\X_\infty |\leq \e$. By Proposition \ref{th3},
we have as $\de \rightarrow 0^+$
 \begin{align} \label{eq4.42}
\max_{|z|\leq \de, |\eta|=1}&
\left(\xi \vee \left[\frac{u(z+x)-u(x)-\P z-\frac{1}{2}\X_{\infty}
:z\ot z}{|z|^2}\right]:\eta
\ot \eta\right) \nonumber\\
& \leq \e \, + \, \max_{|z|\leq \de, |\eta|=1}
\left(\xi \vee\left[ \frac{u(z+x)-u(x)-\P z - \frac{1}{2}\X_{m(\e)}:
z\ot z}{|z|^2}\right]:\eta \ot \eta \right)\\
 &  \leq \ \e \, + \, o(1). \nonumber
 \end{align}
By passing to the limit $\de \rightarrow
0^+$ in \eqref{eq4.42} and letting $\e \rightarrow 0^+$, we
obtain $(P,\X_\infty) \in J^{2,\xi}u(x)$, as a result of Theorem
\ref{th3}. Finally, for any $\eta \in \R^N$, $z \in \R^n$, we have
 \begin{align} \label{eq4.43}
\big[\xi \vee \big(\X + \xi\ot A \big)\big]&: (\eta \ot z) \ot  (\eta \ot z)  \nonumber\\
 & = \ \big(\xi \vee (\X : z \ot z)\big) : \eta\ot \eta  \, + \,  ( \xi^\top \eta)^2(A :  w \ot w) \\
& \geq \ \big(\xi \vee (\X : z \ot z)\big) : \eta\ot \eta  \nonumber
 \end{align}
Consequently, \eqref{eq4.41} follows.      \qed

\ms

Lemma \ref{le9} at the end of Section \ref{section6} supplements
Proposition \ref{pr2} by showing how we can modify the contact jet
$J^{2,\xi}$ along directions perpendicular to $\xi$, that is, when
we can add to $(\P,\X)$ elements of the form $(0,\eta \ot I)$ for
$\eta\, \bot\, \xi$. Now we give an explicit concrete example of jets.

\begin{example}[Calculation of contact jets, cf.\ \cite{K}] \label{ex2}

Let $u :\R \larrow \R^N$ be given by
 \beq \label{eq2.10}
u(z) \ := \ -Az\, \chi_{(-\infty,0]}(z) \  + \ \left(B z
+\frac{C}{2}z^2\right)\chi_{(0,+\infty)}(z),
 \eeq
where $A,B,C \in \R^N$, $A+B \neq 0$. The contact jets of $u$ at zero
are
 \begin{align} \label{eq2.11}
J^{1,\xi}u(0) \ &= \ \left\{
\begin{array}{l}
\O, \hspace{133pt} \xi \neq \frac{A+B}{|A+B|},\\
\left\{\frac{B-A}{2}+t\frac{B+A}{2} \ : \ t\in[-1,+1] \right\}, \xi
= \frac{A+B}{|A+B|},
\end{array}
\right.
\\
\label{eq2.12}
J^{2,\xi}u(0) \ &= \ \left\{
\begin{array}{l}
\O, \hspace{147pt} \xi \neq \frac{A+B}{|A+B|},\\
\left\{\left(\frac{B-A}{2}+t\frac{B+A}{2},\X \right)\ : \ (t,\X) \in
S \right\}, \xi = \frac{A+B}{|A+B|},
\end{array}
\right.
 \end{align}
where
 \begin{align} \label{eq2.13}
S\ :=& \ \Big((-1,+1) \by \R^N \Big) \bigcup \Big(\{-1\}\by
\{C-s(A+B):s\geq 0\}\Big) \nonumber\\
& \ \bigcup \Big(\{+1\}\by \{-s(A+B):s\geq 0\}\Big).
 \end{align}
The proof of the above facts follows by a simple but lengthy computation by using directly the definition of contact jets.
\end{example}

\section{ Ellipticity and consistency with classical notions}
\label{section3}

Now we introduce the appropriate notion of ellipticity
for fully nonlinear second order PDE systems and establish
compatibility between classical and CS.

 \bd[Degenerate elliptic second order systems] \label{de3.1}

Let $u : \R^n \supseteq \Om \larrow \R^N$ be a $C^2$ map. The PDE system \eqref{3.1} is called \emph{degenerate elliptic} when for all $(x,\eta,\P) \in
\Om \by \R^N \by (\R^N \! \ot\R^n)$ the map $\F(x,\eta,\P,\, . \,
) : \R^N \! \ot\R^{n\by n}_s \larrow \R^N$ is monotone, in the sense that the following matrix inequality holds
 \beq \label{eq3.2}
 \Big(\F\big(x,\eta,\P,\X \big)  -  \F\big(x,\eta,\P,\Y \big)\Big)^\top
 \big(\X - \Y \big) \geq \, 0
 \eeq
for all $\X$, $\Y \in \R^N \! \ot\R^{n\by n}_s$, namely
\[
 \Big(\F\big(x,\eta,\P,\X \big) \, - \, \F\big(x,\eta,\P,\Y \big)\Big)^\top
 \big((\X - \Y ):w \ot w \big) \geq \, 0, \ \ \ \text{ for all }w \in \R^n.
 \]
 \ed
By restricting \eqref{eq3.2} to the cases of $N=1 \leq n$ and of $n=1 \leq N$, we recover standard
monotonicity notions which \eqref{eq3.2} extends to the general case. If $N=1 \leq n$ then \eqref{eq3.2} reduces to the standard ellipticity of VS up to a change of sign depending on the convention (see e.g.\ \cite{CIL})
 \beq \label{eq3.4}
\X \leq \Y \ \ \Larrow \ \ \F\big(x,\eta,\P,\X \big) \,\leq \,
\F\big(x,\eta,\P,\Y \big)
 \eeq
If $n=1 \leq N$, then \eqref{eq3.2} reduces to the standard monotonicity of
maps $\R^N \larrow \R^N$. We now derive a characterisation of
Definition \ref{de3.1} which is the form of ellipticity we will actually employ in our analysis.
 \bl \label{le3}
Let $\G: \R^N \! \ot\R^{n\by n}_s \larrow \R^N$. Then, the following are
equivalent:

 \noi (i) For all $\xi \in \mS^{N-1}$ and all $\X$,
  $\Y \in \R^N \! \ot\R^{n\by n}_s$, we have
 \beq \label{eq3.7}
 \xi \vee (\X \, - \, \Y) \ \leq_\ot 0 \ \ \ \Larrow \ \
 \ \xi^\top \big(\G(\X ) \, - \, \G(\Y
)\big) \, \leq\, 0.
 \eeq

\noi (ii) For all $\X$, $\Y \in \R^N \! \ot\R^{n\by n}_s$, we have
 \beq \label{eq3.6}
 \big(\G(\X ) \, - \, \G(\Y )\big)^\top  \big(\X \, - \, \Y \big)\, \geq \, 0.
 \eeq
 \el

\BPL \ref{le3}. By Proposition \ref{pr1}, \eqref{eq3.7} is
equivalent to
 \beq \label{eq3.8}
 \left.
 \begin{array}{c}
\xi^\bot \big(\X \, - \,Y \big) \, = \, 0,\\
\xi^\top \big(\X \, - \,Y \big)  \, \leq  \, 0 
 \end{array}
 \right\}
 \ \ \Larrow \ \ \xi^\top \big(\G(\X ) \, -
\, \G(\Y )\big)\, \leq\, 0.
 \eeq
Assuming \eqref{eq3.8}, we have $\xi^\top(\X \, - \, \Y) \leq 0$ and $\xi \ot \xi^\top(\X  \, -\, \Y)= \X \,-\, \Y$ and also $\xi^\top \big(\G(\X ) \, - \, \G(\Y )\big) \leq 0$. These relations yield
\begin{align} \label{eq3.10}
0 \ & \leq \ \left(\xi^\top \big(\G(\X )
\, - \, \G(\Y )\big)\right) \left(\xi^\top(\X \, - \, \Y) \right) \nonumber\\
 & = \  \big(\G(\X )
\, - \, \G(\Y )\big)^\top \left( \xi \ot \xi^\top(\X \, - \, \Y) \right) \nonumber\\
 & = \  \big(\G(\X )
\, - \, \G(\Y )\big)^\top \left( \X \, - \, \Y\right).\nonumber
 \end{align}
Hence, we obtain \eqref{eq3.6}. Conversely, assuming \eqref{eq3.6}
and that $\xi \vee (\X -  \Y) \leq_\ot 0$, by
Proposition \ref{pr1} we have $ \X - \Y= \xi \ot \xi^\top(\X  - \Y)$ and  hence we get
 \[
0 \,  \leq \,  \big(\G(\X ) \, - \, \G(\Y )\big)^\top (\X \, - \, \Y ) \, = \, \big(\big(\G(\X) \, - \, \G(\Y)\big)^\top\xi\big) \big(\xi^\top(\X
\, - \, \Y)\big).
 \]
Since $\xi^\top(\X - \Y) \leq 0$, we deduce that
$\xi^\top (\G(\X )-\G(\Y )) \leq 0$, as claimed. 
 \qed

\ms

\noi The main result of this section is that CS and classical solutions are compatible for fully nonlinear second order systems which are degenerate elliptic and have continuous coefficients.

 \bt[Consistency] \label{th1} Let $u : \R^n \supseteq \Om \larrow \R^N$ be a continuous map and consider second order system \eqref{3.1}.
 
\smallskip

\noi (a) If $u$ is a contact solution of \eqref{3.1} and the \emph{nonlinearity $\F$ is continuous}, then $u$ solves \eqref{3.1} classically at points of twice differentiability. 

\smallskip

\noi (b) If $u$ is a twice differentiable solution of \eqref{3.1} and the \emph{nonlinearity $\F$ is degenerate elliptic}, then $u$ is a contact solution of \eqref{3.1}.
 \et

\BPT \ref{th1}. (a) If $u:\R^n \supseteq \Om \larrow \R^N$ is a contact
solution of \eqref{3.1}, then, in view of Theorem \ref{th2}, if
$u$ is twice differentiable at $x \in \Om$ we have that
$(\D u (x),\D^2 u   (x))\in J^{2,\pm \xi}u(x)$ for any $\xi \in \mS^{N-1}$.
Hence, by Lemma \ref{l1a}, we have
\beq \label{eq4.22a}
(\pm\xi)^* \F\big(x,u(x),\D u (x),\D^2 u   (x)\big) \, \geq \, 0.
 \eeq
Since $\F$ is continuous, $\xi$-envelopes coincide with $\xi$-projections and we obtain
 \beq \label{eq4.22}
0 \, \leq\, \xi^\top \F\big(x,u(x),\D u (x),\D^2 u   (x)\big) \, \leq\, 0,
 \eeq
for all $\xi \in \mS^{N-1}$. Since $\xi$ is arbitrary, we deduce
that $u$ solves \eqref{3.1} classically. 

\smallskip

\noi (b) Suppose $u$ is a twice differentiable solution of \eqref{3.1} and $\F$ satisfies  \eqref{eq3.2}. Then, if $(\P,\X) \in J^{2,\pm \xi}u(x)$, by Theorem \ref{th2} we have $\P =\D u (x)$ and moreover $\xi \vee
\big(\D^2 u   (x) -\X\big)\leq_\ot 0$. By applying Lemma
\ref{le2}, it follows that whenever $(\P,\X) \in J^{2,\xi}u(x)$,
 \begin{align} \label{eq4.23}
  0 \ &=  \ \xi^\top \F\big(x,u(x),\D u (x),\D^2 u   (x)\big) \nonumber\\
  &\leq  \ \xi^\top \F\big(x,u(x),\D u (x),\X\big) \\
   & \leq \ \xi^* \F\big(x,u(x),\P,\X\big).  \nonumber
 \end{align}
Thus, \eqref{eq4.23} and Lemma \ref{l1a} imply
that $u$ is a contact solution of \eqref{3.1}. \qed

\smallskip

A particular important class of second order PDE systems to which the theory
applies (and has partly been motivated by) is that of
quasilinear ones in non-divergence form. For 
 \begin{align} \label{eq3.12}
A \, = \, A(x,\eta,\P)\ : \ \Om \by \R^N \by(\R^N
\ot \R^n) & \larrow \R^{Nn \by Nn},\\
B \, = \, B(x,\eta,\P)\ : \ \Om \by \R^N \by(\R^N \! \ot\R^n) & \larrow
\R^{N},
 \end{align}
the general form of such systems is
\beq \label{eq3.15}
A \big(\cdot ,u,\D u \big) : \D^2 u    \,+ \, B \big(
\cdot ,u,\D u \big)\, = \, 0
 \eeq
(cf.\ \eqref{1.A}, \eqref{1.B}). According to the next Lemma, in the quasilinear case of
\eqref{eq3.15} the condition of degenerate ellipticity is equivalent
to the rank-one positivity of $A$. The latter is the (weak) Legendre-Hadamard condition, when $A$ is symmetric, namely when $A\in \R^{Nn \by Nn}_{s}$.

 \bl \label{le4} Suppose that $A \in \R^{Nn \by Nn} $. Then, the
linear map $\X \mapsto A:\X$ from $\R^N \! \ot\R^{n\by n}_s$ to
$\R^N$ is monotone if and only if $A \geq_{\ot}0$.
 \el

We note that symmetry of $A$ is not required for this equivalence.
 
\BPL \ref{le4}. The monotonicity of $A$ reads
$\big(A:\X\big)^\top\X \geq 0$. Let us fix
$\eta \in \R^N$, $w\in \R^n$ and set $\X := \eta \ot w \ot w$. Then,
we have
 \begin{align}
0\ \leq & \ \big(A : (\eta \ot w \ot w)\big)^\top (\eta \ot w
\ot w):w\ot w \nonumber\\
= & \ A_{\al i \be j} :  (\eta_\al w_i) (\eta_\be w_j)
(w_k w_k)^2  \nonumber\\
=& \ |w|^4\, A: (\eta \ot w) \ot  (\eta \ot w).  \nonumber
 \end{align}
Hence, $A \geq_{\ot}0$. Conversely, if $ A : (\eta
\ot w) \ot  (\eta \ot w)\geq 0$ for all $\eta \in \R^N$ and all $w\in
\R^n$, we suppose that $\xi \vee \X \leq_\ot 0$ for
some $\xi \in \R^N$. Then, by Proposition \ref{pr1} we have $\X =
\xi \ot X$ where $X:=\xi^\top \X \leq 0$. If $X^{1/2}$
is the symmetric square root of $X$, then we have that $X_{ij}=
X_{ik}^{1/2}X_{jk}^{1/2}$. Hence, $X$ is a sum of positive matrices
$w^{(k)} \ot w^{(k)}$ with $w^{(k)} := X_{ik}^{1/2}e_i \in \R^n$.
Hence, we have
 \[
\xi^\top \big(A : \X \big) \
               =  \ \xi_\al A_{\al i\be j} \xi_\be X_{ij} 
               =  \  A_{\al i\be j}
              ( \xi_\al  w_i^{(k)}) (\xi_\be w_j^{(k)})
              = \, A:
              ( \xi \ot  w^{(k)}) \ot (\xi \ot w^{(k)})
 \]
and therefore $\xi^\top \big(A : \X \big)\leq 0$. By Lemma \ref{le3}, the map $\X \mapsto A:\X$ is monotone.
 \qed
\ms

Now we construct a large class of fully nonlinear systems
which satisfies the ellipticity condition \eqref{eq3.2}.

\begin{example}[Fully nonlinear degenerate elliptic systems] For any nonlinearity $\F=(\F_1,...,\F_N)^\top$ with components of the form
 \beq
\F_\al\ :\ \ \R^n \by \R^N \by (\R^N \! \ot\R^n) \by \R^n \larrow \R^N
 \eeq
such that $F_\al (x,\eta,\P;\cdot )$ is odd and each $l_j \mapsto \F_\al(x,\eta,\P;l_1,...,l_j,...,l_n)$ is homogeneous and increasing for all indices $j,\al$, the next system is degenerate elliptic:
 \beq \label{5.25}
\F_\al\big(\cdot ,u,\D u ,\si(\D^2 u   _\al)\big) \, = \, 0.
 \eeq
\end{example}  

The claim above follows by the next result.

\bl[Monotone functions of the eigenvalues of the hessian] \label{le3.7} Suppose that $g : \R^n \larrow \R^N$ is odd with each component $g_\al$ homogeneous. Suppose further each function $l_j \mapsto g_\al(l_1,...,l_j,$ $...,l_n)$ is increasing, for all indices $j,\al$. Consider
 \beq
\G\ :\ \ \R^N \! \ot\R^{n\by n}_s\larrow \R^N ,\ \ \ \ \G_\al(\X) := g_\al \big(\la_1(X_\al),...,\la_n(X_\al)\big)
 \eeq
where $X_\al := \X_{\al ij}e_i \ot e_j$ and $\big\{\la_1(X_\al),...,\la_n(X_\al)\big\}$ denotes the eigenvalues of $X_\al$, placed in increasing order. Then, $\G$ is monotone in the sense of \eqref{eq3.7}.
\el

\BPL \ref{le3.7}. Fix $\xi \in \mS^{N-1}$, $w\in \R^n$, $\X,\Y \in \R^N \! \ot\R^{n\by n}_s$, an index $\al$ and suppose that $\xi \vee\big(\X -\Y\big) \leq_\ot 0$. Then, we have
 \begin{align} \label{eq3.31}
 0 \ & \geq \ \xi \vee \big(\X  -  \Y\big):(e_{\hat{\al}} \ot w)\ot(e_{\hat{\al}} \ot
 w) \nonumber\\
     & = \ \frac{1}{2}\Big[\xi_\ga \big(\X -\Y \big)_{\be i j } +
      \ \xi_\be \big( \X - \Y \big)_{\ga i j }\Big]\de_{\ga \hat{\al}}w_i \de_{\be
      \hat{\al}}w_j  \nonumber\\
& = \ \xi_{\hat{\al}}\big( \X - \Y \big)_{\hat{\al} i j}w_i w_j \nonumber\\
& = \ \xi_{\hat{\al}} (X-Y)_{\hat{\al}}: w \ot w, \nonumber
 \end{align}
where $\hat{\al} $ denotes free index (no summation). Hence, we obtain $\xi_{\hat{\al}} X_{\hat{\al}} \leq
\xi_{\hat{\al}} Y_{\hat{\al}}$. Since the $k$-eigenvalue function is odd and homogeneous, we have 
$\xi_{\hat{\al}} \la_k(X_{\hat{\al}}) \leq \xi_{\hat{\al}} \la_k(Y_{\hat{\al}})$, for any $\hat{\al}=1,...,N$ and each $k=1,...,n$. Since each $g_\al$ is increasing in each of its arguments, we get
\beq
g_{\hat{\al}} \Big( \xi_{\hat{\al}}\la_1(X_{\hat{\al}}),...,\xi_{\hat{\al}} \la_n(X_{\hat{\al}}) \Big) 
\, \leq\, 
g_{\hat{\al}} \Big( \xi_{\hat{\al}} \la_1(Y_{\hat{\al}}),...,\xi_{\hat{\al}} \la_n(Y_{\hat{\al}}) \Big).
\eeq
Since $g_\al$ is homogeneous and odd, we obtain
\beq
\xi_{\hat{\al}} g_{\hat{\al}} \Big( \la_1(X_{\hat{\al}}),...,\la_n(X_{\hat{\al}}) \Big)
\, \leq\, 
\xi_{\hat{\al}} g_{\hat{\al}} \Big( \la_1(Y_{\hat{\al}}),...,\la_n(Y_{\hat{\al}}) \Big).
\eeq
By summing with respect to $\hat{\al}$, we obtain
$\xi^\top\big( \G(\X) - \G(\Y) \big) \leq 0$. 
 \qed

\ms

Following \cite{CIL}, we can give numerous explicit fully nonlinear degenerate elliptic examples. In particular, the choices $g_\al(l):=l_n$, $g_\al(l):=l_1^{2p_\al +1}$, $g_\al(l):=l_1...l_n$ and $g_\al(l):=(l_1 +...+l_n)^{2p_\al +1}$ lead for any $p_1 ,...,p_N\geq 0$ to the systems
\begin{align}
\max \si(\D^2 u _\al)\ &= \ h_\al(\cdot , u,\D u )\\
\min \big(  \si(\D^2 u _\al)\big)^{2p_\al +1} \ &= \ h_\al(\cdot , u,\D u )\\
\det (\D^2 u _\al)\ &= \ h_\al(\cdot , u,\D u )\ , \ \ u_\al \text{ convex},\\
|\De u_\al|^{2p_\al} \De u_\al\ &= \ h_\al(\cdot , u,\D u ),
\end{align}
which are fully nonlinear and degenerate elliptic for any first order nonlinearity $h$.

\section{ The finer structure of contact jets}
\label{section6}

In this section we study the structure of contact jets more deeply and
demystify the local structure of maps around the point at which a contact jet exists. The principal results are Theorems \ref{th5}-\ref{th6}, which reformulate the matrix inequality defining jets to an \emph{ordinary} inequality \emph{coupling} the $\xi$-projection $\xi^\top u$ and the length of the projection $|\xi^\bot u|$ on the hyperplane normal to $\xi$. The inequality connects a semi-differentiability condition for $\xi^\top u$ (known from the scalar case) to  a new \emph{partial regularity} condition in \emph{codimension-one} for the perpendicular part $\xi^\bot u$. 

 \bt[Structure of first contact jets] \label{th5} Let $u : \R^n \supseteq \Om \larrow
\R^N$ be continuous. Let also $x\in \Om$, $\xi \in \mS^{N-1}$ and $\P
\in \R^N \! \ot\R^n$. Then, the
 following are equivalent:

 \noi(i) $\P \in J^{1,\xi}u(x)$.

\noi(ii) There exists an increasing $\si \in C^1(0,\infty)$
with $\si(0^+)=0$, such that as $z \rightarrow 0$
 \begin{align} \label{eq6.29}
\xi^\top\Big( u(z+x)-u(x) & -\P z\Big)  \, \leq\, -\
\frac{\Big| \xi^\bot\Big( u(z+x)-u(x)-\P z \Big)
\Big|^2}{\si(|z|)|z|} \, + \, \si(|z|)|z|.
\end{align}
 \et

 \bt[Structure of second contact jets] \label{th6} Let $u : \R^n \supseteq \Om \larrow
\R^N$ be continuous. Let also $x\in \Om$, $\xi \in \mS^{N-1}$ and
 $(\P,\X) \in \R^N \! \ot(\R^n \by \R^{n\by n}_s)$. Then, the
 following are equivalent:

 \noi(i) $(\P,\X) \in J^{2,\xi}u(x)$.

\noi(ii) There exists an increasing $\si \in C^2(0,\infty)$
with $\si(0^+)=0$, such that as $z \rightarrow 0$
 \begin{align} \label{eq6.30}
\xi^\top\Big( u(z+x) &  -u(x)  -\P z-\frac{1}{2}\X:z\ot z
\Big)  \nonumber\\
 & \leq  \ -\
\frac{\Big| \xi^\bot\Big(u(z+x)-u(x)-\P z - \dfrac{1}{2}
\X:z\ot z\Big)
\Big|^2}{\si(|z|)|z|^2}\  + \ {\si(|z|)|z|^2}.
\end{align}
 \et

By \eqref{eq6.29} and \eqref{eq6.30} we obtain that the existence of nontrivial contact jets implies a local structure for the map: the codimension-one projection of $u$ on the hyperplane $\xi^\bot$ must be more regular than the projection $\xi^\top u$. Actually there is a bootstrap of regularity between $\xi^\top u $ and $\xi^\bot u$, which balances at differentiability:

 \begin{corollary}[Codimension-one bootstrap regularity imposed by jets]
 \label{cor3} Suppose that $u : \R^n \supseteq \Om \larrow \R^N$ is a continuous
 map, $x\in \Om$, $\xi \in \mS^{N-1}$. 

\smallskip

\noi (1) Let $\P \in J^{1,\xi}u(x)$ and $L(z):=u(z+x)-u(x)-\P z$. Then
\[
\ \ \ \ \ \ \xi^\top L(z)\ =\ O(|z|^\be) \ \ \ \Longrightarrow \ \ \ \xi^\bot L(z)\ =\ o(|z|^{\be+\frac{1-\be}{2}}),
\]
as $z\to 0$,  for any $\be \in[0,1]$. In particular, for $\be \in\{0,1\}$ the following holds: since $\xi^\top u$ is $C^0$ at $x$, $\xi^\bot u$ is $C^{\frac{1}{2}+}$ near $x$. If $\xi^\top u$ is differentiable at $x$, then so is $\xi^\bot u$.

\smallskip

\noi (2) Let $(\P,\X) \in J^{2,\xi}u(x)$ and $Q(z):=u(z+x)-u(x)-\P z-\frac{1}{2}\X:z\ot z$. Then
\[
\ \ \ \  \ \ \xi^\top Q(z)\ =\ O(|z|^\ga) \  \ \ \Longrightarrow \ \ \ \xi^\bot Q(z)\, = \, o(|z|^{\ga + \frac{2-\ga}{2}}),
\]
as $z \to 0$,for any $\ga \in[0,2]$. In particular,  for $\ga\in \{0,2\}$ the following hold: since $\xi^\top u$ is $C^0$ at $x$, $\xi^\bot u$ is $C^{1}$ near $x$. If $\xi^\top u$ is twice differentiable at $x$, so is $\xi^\bot u$.
 \end{corollary}

\BPCOR \ref{cor3}. We rewrite \eqref{eq6.29} and \eqref{eq6.30} as
 \begin{align} \label{eq6.30a}
 \big|\xi^\bot L(z) \big|^2\ &\leq   \ \si(|z|)|z| \Big(-\xi^\top L(z)\, + \, \si(|z|)|z|\Big)\\
 \big|\xi^\bot Q(z) \big|^2 \ &\leq   \ \si(|z|)|z|^2 \Big(-\xi^\top Q(z)\, + \, \si(|z|)|z|^2\Big)
\end{align}
and the desired conclusions readily follow.   \qed

\medskip

The fact of existence of nowhere improvable H\"older functions implies

 \begin{corollary} \label{cor4} For any $\xi \in \mS^{N-1}$, there exists a map $u \in (C^0
 \set C^{\frac{1}{2}})(\R^n)^N$
such that $u$ does not possess nontrivial first $\xi$-jets anywhere. Similarly, for any $\xi \in \mS^{N-1}$, there exists a map $u \in (C^0\set C^{0,1})(\R^n)^N$ such that $\big\{ x\in \R^n \ | \ J^{2,\xi}u(x)\neq \emptyset \big\}  =  \emptyset$.
 \end{corollary}
Hence, obstructions arising in the vectorial case imply that first contact jets are efficient for H\"older $C^{\frac{1}{2}}$ maps and second contact jets are efficient for Lipschitz $C^{0,1}$ maps. In the scalar case obstructions disappear and semi-jets $J^{2,\pm}$ are efficient for merely $C^0$ functions. We interpret this fact by saying that ``in the vectorial case only $1/2$ of the derivatives can be interpreted weakly, the rest $1/2$ must exist classically".

In order to prove Theorems \ref{th5}-\ref{th6}, we need a technical tool. Let $R \in \R^N$ and $\xi\in \mS^{N-1}$. By Lemma \ref{le1}, the tensor product $\xi \vee R$ is a rank-two symmetric tensor.

  \bl[Representations of the spectrum of symmetrised tensor products]
\label{le7} For any $R\in \R^N$, set $s(R):=\ 2\big(\sgn(\xi^\top R)\big)^+ -1$. Then, we have the identities
 \begin{align} \label{eq6.2}
 \max \si (\xi \vee R) \, &= \, \big(\xi^\top R\big)^+  \ +\
 \frac{|R|}{4}\Big|\sgn(R)-s(R)\xi\Big|^2,
\\
\label{eq6.7}
 \max \si (\xi \vee R) \, &= \, \max \big\{\xi^\top R,0 \big\}  \ +\
 \frac{|R|}{4} \min \Big| \sgn(R)\pm \xi\Big|^2.
 \end{align}
 \el

\BPL \ref{le7}. By observing that for any $a \in \R$, we have
 \begin{align} \label{eq6.4}
2\big(\sgn(a)\big)^+  - 1 \, = \, \big(\chi_{(0,\infty)} -
\chi_{(-\infty,0]}\big)(a),\nonumber
\end{align}
we obtain that $s(R)=\big(\chi_{(0,\infty)} - \chi_{(-\infty,0)}\big)(\xi^\top R)$. We assume $R\neq0$, since \eqref{eq6.2} is trivial if $R=0$. Let
$l(R)$ denote the right hand side of \eqref{eq6.2}. Then, we compute
\[
\begin{split}
l(R)\, & = \, \big(\xi^\top R\big)^+  \ +\
 \frac{|R|}{4}\left|\frac{R}{|R|}- \xi \big(\chi_{(0,\infty)} \, - \,
\chi_{(-\infty,0]}\big)(\xi^\top R) \right|^2    
\\
&= \, \frac{\big|\xi^\top R\big| +
\xi^\top R }{2} \, + \, \frac{|R|}{2}\bigg(1 \, -\,  \frac{\xi^\top
R}{|R|}\big(\chi_{(0,\infty)} \, - \, \chi_{(-\infty,0]}\big)(\xi^\top
R) \bigg)
\\
& = \ \frac{1}{2}\Big[ \big|\xi^\top R\big| + \xi^\top R  \ +
\ |R| \ - \xi^\top R \big(\chi_{(0,\infty)}
\, + \, \chi_{(-\infty,0]}\big)(\xi^\top R) \Big] 
 \\
  & = \ \frac{1}{2}\Big[ \big|\xi^\top R\big| \, + \, |R|
   \ + 2 \xi^\top R \chi_{(-\infty,0]}(\xi^\top R) \Big]
\\
& = \,
\dfrac{1}{2}\Big[ \big|\xi^\top R\big| \, + \, |R|
   \ + 2 \xi^\top R \Big] \chi_{\{\xi^\top R \leq 0\}}(R)
   \,+\,
  \dfrac{1}{2}\Big[ \big|\xi^\top R\big| \, + \, |R|\Big] \chi_{\{\xi^\top R > 0\}}(R),   
 \end{split}
 \]
which gives $l(R)=\max \si(\xi \vee R)$, as a result of Lemma \ref{le1}. This
establishes \eqref{eq6.2}. Let us now establish \eqref{eq6.7}. The
elementary identity
 \beq \label{eq6.8}
\frac{|R|}{4} \Big| \sgn(R) \,-\, \xi\Big|^2 \ =  \ \frac{|R|}{4} \Big| \sgn(R) \,+\, \xi\Big|^2 \, + \, \xi^\top R
 \eeq
implies that
 \begin{align} \label{eq6.9}
 \frac{|R|}{4} \min \Big| \sgn(R)\pm \xi\Big|^2 \ & = \
 \left\{
\begin{array}{l}
\dfrac{|R|}{4} \Big| \sgn(R) \,+\, \xi\Big|^2, \ \ \text{ if } \xi^\top R \leq 0,
 \nonumber \medskip\\
\dfrac{|R|}{4} \Big| \sgn(R) \,-\, \xi\Big|^2 ,
 \ \ \text{ if } \xi^\top R > 0,  \nonumber
\end{array}
\right. \\
&= \ \frac{|R|}{4}\left|\frac{R}{|R|}\, -\, \xi \big(\chi_{(0,\infty)} \,
- \, \chi_{(-\infty,0]}\big)(\xi^\top R) \right|^2 \\
& = \ \frac{|R|}{4}\Big|\sgn(R)\, -\, s(R)\xi\Big|^2. \nonumber
\end{align}
By comparing \eqref{eq6.2} and \eqref{eq6.9}, we see that
\eqref{eq6.7} follows.           
 \qed   
 
 \medskip

The following is the first step towards Theorems \ref{th5}, \ref{th6}.

 \bt[Equivalent formulations of contact jets] \label{th4}
Let $R : \R^n \supseteq \overline{\Om} \larrow \R^N$ be continuous, $0\in \overline{\Om}$, $R(0)=0$, $\xi \in \mS^{N-1}$. Let also $p \in \{1,2\}$. Then, the following statements are
equivalent:

\noi (i) $0 \in J^{p,\xi}R(0)$, that is, $\max\si(\xi\vee R(z))\leq
o(|z|^p)$ as $\overline{\Om} \ni z\rightarrow 0$.

\noi (ii) We have
 \beq \label{eq6.11}
  \xi^\top R (z) \leq o(|z|^p) \ \text{ and } \ \dfrac{\big|\xi^\bot R (z)\big|^2}{|R(z)|} =  o(|z|^p) \ \text{ as
}\overline{\Om} \ni z \rightarrow 0.
 \eeq

\noi (iii) There exist maps $\rho : \R^n \larrow \xi^\bot \sub \R^N$ and $\si : \R^n\larrow [0,\infty)$ satisfying that $\rho(z)=o(|z|^{p/2})$ and $\si(z)=o(|z|^{p})$ as $\overline{\Om} \ni z\rightarrow 0$ and also
  \beq \label{eq6.12}
 \xi^\top R  \leq  \si \ \text{ and } \  \xi^\bot R  = \rho \Big(\frac{1}{2}|\rho|^2 +
\Big(\frac{1}{2}|\rho|^2\Big)^2 + \big|\xi^\top R\big|^2 \Big)^{1/2} \ \ 
 \text{ on } \overline{\Om}.
 \eeq

\noi (iv) If we set $T := \big\{|\xi^\bot R| \leq |\xi^\top R|
\big\} \sub \overline{\Om}$, then
  \beq \label{eq6.13}
  \left\{
\begin{array}{l}
\left. \begin{array}{l}
\ \ \ \ \xi^\top R (z) \, \leq\, o(|z|^p),   \\
\dfrac{\big|\xi^\bot R (z)\big|^2}{\big|\xi^\top R (z)\big|} \ = \
o(|z|^p),
\end{array}\right\} \  \text{ as }\overline{\Om} \cap T \ni z \rightarrow 0, \medskip\\
 \hspace{25pt} |R(z)| \, = \, o(|z|^p),  \hspace{20pt} \text{ as }\overline{\Om} \set T \ni z
\rightarrow 0.
\end{array}
  \right.
 \eeq
  \et

We observe that when $N=1$, then $\xi^\bot \equiv 0$ and we recover a single
inequality along $\xi \ot \xi \cong \R$ which coincides with that
of scalar semijets.

\BPT \ref{th4}. We begin by proving that (i) is equivalent to (ii).
If we assume (i), then by the representation formula \eqref{eq6.2},
it is equivalent to
  \beq \label{eq6.14}
\xi^\top R (z) \leq  o(|z|^p) \ \text{ and } \
|R(z)|\big|\sgn(R(z))\,-\, s(R) \xi \big|^2 = o(|z|^p),
 \eeq
as $\overline{\Om} \ni z \rightarrow 0$. By \eqref{eq6.14}, we have $\big|  R - |R|s(R)\xi\big|(z) \leq o\big(|z|^p |R(z)|^{\frac{1}{2}}\big)$, as $\overline{\Om} \ni z\rightarrow 0$. Hence, there exists an $r: \overline{\Om} \larrow \R^N$ with $|r(0^+)|=0$ such that
 \beq \label{eq6.16}
\big( R- |R|s(R)\xi \big)(z) \, = \, r(z)|z|^{\frac{p}{2}}|R(z)|^{\frac{1}{2}},
 \eeq
on $\overline{\Om}$. By projecting \eqref{eq6.16} onto $\xi^\bot$, we obtain $\xi^\bot R(z)\, = \, (\xi^\bot r)(z)|z|^{\frac{p}{2}}|R(z)|^{\frac{1}{2}}$ on $\overline{\Om}$. We set $\rho(z):= (\xi^\bot r)(z)|z|^{\frac{p}{2}}$. Therefore
 \begin{align} \label{eq6.18}
 \frac{\big|\xi^\bot R(z)\big|^2}{|R(z)|}\ & = \ |\rho(z)|^2\, = \, o(|z|^p), 
\end{align}
as $\overline{\Om} \ni z \rightarrow 0$. By \eqref{eq6.14} and \eqref{eq6.18},
(ii) follows. Conversely, assume (ii). It suffices to verify that $\max \si(\xi\vee
R(z))\leq o(|z|^p)$ as $\overline{\Om} \cap \{R\neq 0\} \ni z\rightarrow 0$. Then, we calculate
\begin{align} 
 \frac{|R|}{4}\Big|\sgn(R)- s(R)\xi\Big|^2 \ & =\ \frac{1}{4|R|}\Big|R- s(R)|R|\xi\Big|^2
\nonumber\\
 & =\ \left\{
 \begin{array}{l}
\dfrac{1}{4|R|}\Big|R\,-\, |R|\xi\Big|^2 ,\text{ on }\overline{\Om}\cap\{\xi^\top
R>0\}\cap \{R\neq0\},
\\
\dfrac{1}{4|R|}\Big|R\,-\,  |R|\xi\Big|^2 ,\text{ on }\overline{\Om}\cap\{\xi^\top R
\leq 0\}\cap \{R\neq0\}. \medskip
\end{array}
 \right.
\\
 & = \left\{
 \begin{array}{l}
\dfrac{\big(|R|\,-\,  \xi^\top R\big)}{2} ,\text{ on }\overline{\Om}\cap\{\xi^\top
R>0\},
\medskip \nonumber
\\
\dfrac{\big(|R| \,+\,  \xi^\top R\big)}{2} ,\text{ on }\overline{\Om}\cap\{\xi^\top R
\leq 0\}\cap \{R\neq0\}. 
\end{array}
 \right.
\end{align}
Since $\xi^\bot =I-\xi \ot \xi$, we have that
\[
\Big(|R(z)|- \xi^\top R(z)\Big)\Big(|R(z)|+ \xi^\top R(z)\Big) = |\xi^\bot R(z)|^2 = o\big(|z|^p|R(z)|\big),  
 \]
as $\overline{\Om} \ni z \rightarrow 0$. Since on $\overline{\Om}\cap\{\xi^\top R>0\}$ we
have $1\leq \big(|R| + \xi^\top R\big)/|R|$ and on
$\overline{\Om}\cap\{\xi^\top R\leq 0\}\cap\{R\neq0\}$ we have $1\leq \big(|R|
- \xi^\top R\big)/|R|$, we infer that
\begin{align} \label{eq6.21}
 \frac{|R|}{4}\Big|\sgn(R)- s(R)\xi\Big|^2 \  \leq \
 \frac{\left(|R|^2\, -\, |\xi^\top R|^2\right)}{2|R|}(z)\, = \, o(|z|^p), 
 \end{align}
as $\overline{\Om} \ni z \rightarrow 0$. Thus, \eqref{eq6.11} and
\eqref{eq6.21} imply \eqref{eq6.14}, which is equivalent to (i). Let us now prove the equivalence between (ii) and (iii). If we
assume (ii), we define
 \beq \label{eq6.22}
\si := \big(\xi^\top R\big)^+ ,\ \ \ \  
\rho :=  \big({\xi^\bot
R}{|R|^{-\frac{1}{2}}}\big)\chi_{\{R\neq0\}\cap\overline{\Om}}.
 \eeq
It follows that $\si$, $\rho$ have the desired properties and by
\eqref{eq6.22}, we have $|R||\rho|^2 = |\xi^\bot R|^2 =|R|^2 - |\xi^\top R|^2$. It follows that $|R|$ is the positive solution of the quadratic equation
$t^2 - |\rho|^2 t -|\xi^\top R|^2=0$. Hence,
 \beq  \label{eq6.24}
|R|\, = \, \frac{1}{2}|\rho|^2 \, + \, \Big(\Big(\frac{1}{2}|\rho|^2
\Big)^2 \, + \, |\xi^\top R|^2\Big)^{1/2}.
 \eeq
Thus, \eqref{eq6.22} and \eqref{eq6.24} imply \eqref{eq6.12}.
Conversely, if we assume (iii) and let $S$ be defined by the formula giving $R$ above, 
then $S$ solves the equation $t^2 - |\rho|^2 t -|\xi^\top
R|^2=0$. Hence, by the above and perpendicularity, we
have
 \beq
S^2\  =\ |\xi^\top R|^2 \, + \, |\rho|^2 S \  =\ |R|^2 \, - \, |\xi^\bot R|^2 \, + \, |\rho|^2 S \, = \, |R|^2. 
\eeq
This, $S= |R|$ and hence $\xi^\bot R=\rho|R|^{\frac{1}{2}}$. As a result, \eqref{eq6.12} implies \eqref{eq6.11} as claimed. We conclude by proving that (iv) is equivalent to (ii). For, let us split $\R^N$ as $\xi^\top \oplus \xi^\bot$ and equip it with the norm $\|R\|:=\max\big\{|\xi^\top R|,|\xi^\bot
R|\big\}$. Then, we have
\begin{align} 
\frac{|\xi^\bot R|^2}{\|R\|}\ &= \ \frac{|\xi^\bot
R|^2}{\max \big\{|\xi^\top R|,|\xi^\bot R|\big\}} = \ \frac{|\xi^\bot R|^2}{|\xi^\top R|} \chi_{ \overline{\Om} \cap T} \ + \
|\xi^\bot R| \chi_{ \overline{\Om} \set T}. 
\end{align}
By the above and norm equivalence on $\R^N$, \eqref{eq6.11} is equivalent to
 \beq \label{eq6.28}
\left\{
\begin{array}{c}
\ \xi^\top R(z) \, \leq\, o(|z|^p), \  \text{ as }\overline{\Om} \ni z
\rightarrow
0, \medskip\\
\left\{
\begin{array}{l}
\dfrac{|\xi^\bot R(z)|^2}{|\xi^\top R(z)|}\, = \, o(|z|^p), \ \text{
as }\overline{\Om}\cap T
\ni z \rightarrow 0, \medskip\\
\ |\xi^\bot R(z)|^2 \, = \, o(|z|^p), \ \text{ as }\overline{\Om}\set T \ni z
\rightarrow 0,
\end{array} \right.
\end{array} \right.
 \eeq
Since on $\overline{\Om} \set T$ we have $|\xi^\top R(z)| < |\xi^\bot
R(z)|=o(|z|^p)$ as $\overline{\Om}\set T \ni z \rightarrow 0$, \eqref{eq6.28}
is equivalent to \eqref{eq6.13} and the theorem follows.      \qed
\medskip

\BPT\!\!\textbf{s} \ref{th5}, \ref{th6}. By Theorem \ref{th4}, it suffices to prove the following:
 \bc \label{cl1} If $R:\R^n \larrow \R^N$ is continuous,
 $R(0)=0$, $\xi \in \mS^{N-1}$,  $p \in \{1,2\}$, then
 \beq \label{eq6.33}
  0 \in J^{p,\xi}R(0) \ \ \ \
\Longleftrightarrow \ \ \ \ \xi^\top R(w) - (\si(|w|)|w|)^p \, \leq\, -
\ \frac{\big|\xi^\bot R(w)\big|^2}{(\si(|w|)|w|)^p}
 \eeq
as $w\rightarrow 0$, for some  increasing $\si \in C^p(0,\infty)$
which satisfies $\si(0^+)=0$.
 \ec

\noi To this end, assume $0\in J^{p,\xi}R(0)$. Then, by Theorem
\ref{th3}, there exists an  increasing $\rho \in C^p(0,\infty)$ with $0\leq
\rho(w)=o(|w|^p)$ as $w\rightarrow 0$ such that
 \beq \label{eq6.34}
\xi \vee R(w) \, \leq\, \rho(w)I.
 \eeq
Let $\{\xi_1,..,\xi_{N-1}\}$ be an
orthonormal base of the hyperplane $\xi^\bot \sub
\R^N$. Then, $R$ and the identity $I$ can be written as:
\beq
R\, = \, (\xi^\top R)\xi \, + \sum_{\al=1}^{N-1}({\xi_\al}^\top R)
\xi_\al ,\ \ \ 
I\, = \, \xi \ot \xi \, + \sum_{\al=1}^{N-1}\xi_\al \ot \xi_\al .  \label{eq6.35}
\eeq
By plugging \eqref{eq6.35} into \eqref{eq6.34},
we obtain
 \beq  \label{eq6.37}
\xi \vee \left[ (\xi^\top R)\xi \, +\, \sum_{\al=1}^{N-1}({\xi_\al}^\top R)
\xi_\al \right] \, \leq\, \rho \, \xi \ot \xi \ +\ \rho
\sum_{\al=1}^{N-1}\xi_\al \ot \xi_\al.
 \eeq
By applying ``$:\xi \ot \xi$'' to \eqref{eq6.37} and employing orthonormality of the base, we infer that $-\ \xi^\top R  +  \rho \geq  0.$. Let now $t\in\R\set\{0\}$ and $\be \in \{1,...,N-1\}$ be
 fixed and apply again ``$:(t\xi_\be +\xi) \ot (t\xi_\be +\xi)$'' to
\eqref{eq6.37} to obtain
 \begin{align}  \label{eq6.39}
|\xi|^4 \xi^\top R \, + \sum_{\al=1}^{N-1} \big(\xi_\al^\top R \big)
 \Big(\xi_\al^\top(t\xi_\be +\xi) \Big)
 \Big(\xi^\top(t\xi_\be +\xi) \Big)  \leq \, \rho|\xi|^4  + \, \rho \sum_{\al=1}^{N-1}
\big(\xi_\al^\top(t\xi_\be +\xi)\big)^2. \nonumber
 \end{align}
By orthogonality of the base, we deduce that $t \xi_\be^\top R  \leq \rho t^2 +  \rho  +  \xi^\top R$. Since this holds for both $\pm t$, we infer that
 \beq  \label{eq6.41}
\big|\xi_\be^\top R\big| \, \leq\, \rho |t| \, + \, \frac{\rho \ - \
\xi^\top R}{|t|}
 \eeq
and the choice $t:= \left(({\rho - \xi^\top
R}{)/\rho}\right)^{{1}/{2}}$ in \eqref{eq6.41} implies $\big|\xi_\be^\top R\big|^2  \leq 4 \rho \big(\rho + \xi^\top R\big)$. By summing with respect to $\be$, we obtain
\begin{align} \label{eq6.43}
\big|\xi^\bot R\big|^2 & = \, \sum_{\be=1}^{N-1} \big|\xi_\be^\top R\big|^2 \leq \, 4 (N-1)  \Big( - \xi^\top R \, + \, 4 (N-1)\rho \Big)\rho .
\nonumber
\end{align}
The above estimate implies the direction ``$\Longrightarrow$'' of
Claim \ref{cl1} for the choice $\si(|w|):={|w|}^{-1}(4(N-1)\rho(|w|))^{1/p}$. Conversely, assume the validity of the inequality in \eqref{eq6.33}
for such a $\si$ and set $\rho(w):=(\si(|w|)|w|)^p$. Then, we have
\beq \label{eq6.44}
 \big|\xi^\bot R\big|^2  \leq\, \rho \big(- 
\xi^\top R \, + \, \rho\big),
 \eeq
locally near $0\in \R^n$. Since $\rho>0$ near zero, \eqref{eq6.44}
readily gives $\xi^\top R \leq  \rho(w) = o(|w|^p)$ as $w\rightarrow 0$. By setting
$T  := \big\{\big|\xi^\bot R\big|\leq \big|\xi^\top R\big|\big\}$ and
$\Om  := \big\{\big|\xi^\top R\big|>\rho\big\} $, the inequality \eqref{eq6.44} implies on $T \cap \Om$ that
 \begin{align}  \label{eq6.47}
\frac{\big|\xi^\bot R\big|^2}{\big|\xi^\top R\big|}(w)\ & \leq \
\frac{\rho(w)\big(- \ \xi^\top R(w) \, + \, \rho(w)\big)}{\big|\xi^\top R(w)\big|}
  \leq \ \rho(w)\, + \,  \frac{\rho^2(w)}{\big|\xi^\top R(w)\big|} \leq \ 2 \rho(w),  \nonumber
\end{align}
as $T\cap \Om \ni w\rightarrow 0$. Hence, by the implication (iv)
$\Longrightarrow$ (i) of Theorem \ref{th4}, we obtain
 \beq \label{eq6.48}
\max \si\big(\xi \vee R(w)\big)\, \leq\, o(|w|^p),
 \eeq
as $T \cap \Om \ni  w\rightarrow 0$. On the other hand, on $T \set \Om$ we have $\big|\xi^\top R\big|\leq \rho$ and also $\big|\xi^\bot R\big| \leq \big|\xi^\top R\big|$, hence by Lemma \ref{le1} we estimate
 \begin{align} \label{eq6.49}
\max \si\big(\xi \vee R(w)\big)\ & = \ \frac{1}{2}\Big(|R(w)| \, + \,
\xi^\top R(w)\Big) \leq \, |R(w)| \,=  \nonumber\\
& = \ \Big(\big|\xi^\top R\big|^2+ \, \big|\xi^\bot
R\big|^2\Big)^{\frac{1}{2}} \leq \, \big(\rho^2(w) + \rho^2(w) \big)^{\frac{1}{2}}\,  = \, \sqrt{2} \rho(w) ,  \nonumber
\end{align}
as $T\set \Om \ni w \rightarrow 0$. Therefore, we deduce that
 \beq \label{eq6.50}
\max \si\big(\xi \vee R(w)\big)\, \leq\, o(|w|^p),
 \eeq
as $T \ni w\rightarrow 0$. Now, by
\eqref{eq6.44} on $\R^n \set T =
\big\{\big|\xi^\top R\big|<\big|\xi^\bot R\big|\big\}$  we have that
 \begin{align} \label{eq6.51}
\big|\xi^\bot R\big|^2 & \leq \ \rho\big(- \xi^\top R\, + \, \rho\big)   \, \leq\, \rho \big|\xi^\top R \big|\, + \, \rho^2  \, \leq\, \rho \big|\xi^\bot R \big|\, + \, \rho^2. 
\end{align}
Hence, it holds that $\big|\xi^\bot R\big|^2 -  \rho \big|\xi^\bot R \big| -  \rho^2 \leq 0$
 and by comparing $\big|\xi^\bot R\big|$ with the solutions of the
binomial equation $t^2 - \rho t - \rho^2 = 0$, we find
 \begin{align} \label{eq6.53}
\big|\xi^\bot R(w)\big| & \leq \ \frac{1+\sqrt{5}}{2}\rho(w) \, = \, o(|w|^p), 
\end{align}
as $\R^n \set T \ni w\rightarrow 0$. Thus, by employing again Lemma
\ref{le1}, we estimate
 \beq
 \begin{split} \label{eq6.54}
\max \si\big(\xi \vee R(w)\big) & =  \frac{1}{2}\Big(|R(w)| +
\xi^\top R(w)\Big) \leq \, \Big(\big|\xi^\top R\big|^2 +\big|\xi^\bot
R\big|^2\Big)^{\frac{1}{2}} \leq \frac{1+\sqrt{5}}{\sqrt{2}}\rho(w),  
\end{split}
\eeq
as $\R^n \set T \ni w \rightarrow 0$. By estimates \eqref{eq6.50}
and \eqref{eq6.54} we deduce $\max \si\big(\xi \vee R(w)\big)
\leq o(|w|^p)$ as $w \rightarrow 0$ and consequently $0\in J^{p,\xi}R(0)$. As a result, Claim \ref{cl1} follows. \qed

\medskip

The following result certifies that the projections along
$\xi^\bot$ of second contact $\xi$-jets are ``stiff'' and if the map is twice differentiable, no variations can be performed.

 \bl \label{le9} Let $u :\R^n \subseteq \Om \larrow \R^N$ be continuous and fix $x\in \R^n$, $\xi \in \mS^{N-1}$, $\eta \in \mS^{N-1}\cap
 \xi^\bot$ and $A \in \R^{n\by n}_+$. Then, if
 $(\P,\X) \in J^{2,\xi}u(x)$, we have
  \beq
\left(\eta^\top\P\, ,\, \eta^\top\X-\frac{1}{2}A\right) \in
J^{2,+}(\eta^\top u)(x)\ \ \Longrightarrow \ \ (\P,\X -\eta \ot A)
\in J^{2,\xi}u(x).
   \eeq
 \el

Lemma \ref{le9} says that we can modify $(\P,\X)$ by adding an element of the form $(0,-\eta \ot A)$ with $A\geq 0$ along a direction $\eta$ in the normal hyperplane $\xi^\bot$ if we can add the element $\left(0,-\frac{1}{2}A\right)$ from the superjet $J^{2,+}$ of the projection $\eta^\top u$. For modifications along the $\xi$-direction, see Proposition \ref{pr2}.

\BPL \ref{le9}. We set
 \beq
Q_{\P,X}(z)\ := \ u(z+x) - u(x)- \P z- \frac{1}{2}\X:z
\ot z.
 \eeq
By employing that $\eta\, \bot\, \xi$, we estimate
\begin{align}
\Big|\xi^\bot \Big( Q_{\P,X}(z)  -\eta \ot
A:z\ot z\Big)\Big|^2  =& \ \left|\xi^\bot \big(u- Q_{\P,X}\big)(z)\right|^2
 \, + \, \frac{1}{4}\big|A:z \ot z\big|^2 \nonumber\\
& + \ \big(A: z \ot z\big) \, \eta^\top \big(u- Q_{\P,X}\big)(z)  \nonumber \\
\leq & \ o(|z|^2)\Big[-\xi^\top \big(u- Q_{\P,X}\big)(z) \, + \, o(|z|^2)\Big]\\
& + \ \big(A:z \ot z\big)\Big[ \, \frac{A}{4}:z
\ot z+ \eta^\top \big(u- Q_{\P,X}\big)(z)\Big],\nonumber
\end{align}
as $z\to 0$. Hence, by assumption we have
\begin{align}
\Big|\xi^\bot \Big( Q_{\P,X}(z) & -\eta \ot
A:z\ot z\Big)\Big|^2 \, \leq\, o(|z|^2)\Big[-\xi^\top \big(u- Q_{\P,X}\big)(z) \, + \, o(|z|^2)\Big] \nonumber\\
& + \ \big\|A\big\||z|^2\eta^\top \Big( u(z+x) -u(x)-\P z \ - \frac{1}{2}\left(\X - \frac{\eta}{2} \ot A \right):z \ot z\Big)\\
\leq &\ o(|z|^2)\Big[-\xi^\top \big(u- Q_{\P,X}\big)(z) \ + \
o(|z|^2)\Big] \, + \, o(|z|^4), \nonumber
\end{align}
as $z\to 0$. By increasing the $o(1)$ functions appearing in the
summands appropriately, we incorporate the $o(|z|^4)$ term in the
first summand and therefore obtain that $(\P,\X -\eta \ot A) \in
J^{2,\xi}u(x)$, as claimed. \qed

\ms

\section{ The extremality notion of contact maps}
\label{section7}

So far, our central objects of study have been contact jets, a certain type of generalised pointwise derivatives. Jets in fact introduce in an implicit non-trivial fashion an extremality notion for maps, which we will now exploit. This notion extends  min and max of scalar functions to the vector-valued case, effectively extending the ``Maximum Principle calculus" ($\D u =0$ and $\D^2 u   \leq0$ at maxima of $u$) to the vectorial case. This device allows the ``nonlinear passage of derivatives to test maps". This extremality notion, although simple in its form, presents peculiarities and is not obvious how it arises. Hence, we have chosen to base the PDE theory of CS to jets rather than to extrema, since jets seem more reasonable due to the formal resemblance to their scalar counterparts.

\smallskip

\noi {\bf Motivation}. We begin by motivating the notions that follow. Let $u : \R \larrow \R^N$ be a smooth curve. Every reasonable definition of extremal point $u(\bar{x}) \in \R^N$ at $\bar{x}\in\R$ must imply that $|u'(\bar{x})|=0$. However, this is impossible if $N\geq 2$ as the example of unit speed curves certifies for which $|u'|\equiv 1$. In order to succeed we must radically change our point of view of ``extremals".  The idea is to relax the pointwise notion to a flexible \emph{functional notion of ``extremal map''} which takes into account the possible ``twist''. Our viewpoint is the following: if $N=1$ and $u : \R \larrow \R$ has a maximum $u(\bar{x}) \in \R$ at $\bar{x}\in\R$, then we can identify the extremum $u(\bar{x})$ with the constant function $\psi \equiv u(\bar{x}) : \R \larrow \R$ which passes through $\bar{x}$  (Figure 1(a)).
\[
\underset{\text{Figure 1(a). \hspace{100pt} Figure 1(b).}}{\includegraphics[scale=0.4]{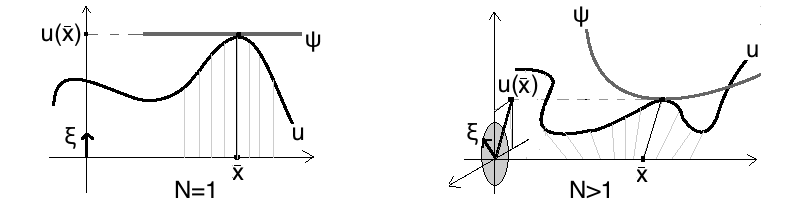}}
\]
When $N\geq 2$ we can view extrema as maps $\psi : \R \larrow \R^N$ passing at $x$ through $u(\bar{x})$ which generally are \emph{nonconstant} (Figure 1(b)).  

Going back to $N=1$, we see that maximum can be viewed as a constant function $\psi$ ``touching $u$" at $x$ in the direction $\xi=+1$ and minimum as ``touching $u$" at $x$ in the direction $\xi=-1$. When $N\geq 2$, there still exists a stiffer notion of ``touching $u$" at $x$ by a map $\psi$ in a unit direction $\xi \in \R^N$. We will call this stronger touching notion \emph{contact}. 

There are two intriguing properties associated with contact which are source of difficulties. Firstly, the notion of contact comprises \emph{a notion of extremum not connected to any ordering of $\R^N$} (when $N\geq 2$). Its main utility is the ``nonlinear passage of derivatives to test maps" in our PDE systems. Secondly, the contact \emph{has order} itself: roughly, ``first order contact" implies ``gradient equality" and ``second order contact" implies an appropriate ``hessian inequality".

To begin, let $C_x$ denote a generic cone function with vertex at $x\in \R^n$ and some slope $L>0$, that is
$C_x(z):=L|z-x|$.

\bd[Contact maps] \label{de7.3} Let $u : \R^n \supseteq \overline{\Om} \larrow
\R^N$ be continuous and fix $x\in \overline{\Om} $ and $\xi \in \mS^{N-1}$. 

\smallskip

\noi (1) \emph{The map $\psi \in C^1(\R^n)^N$ is a first contact $\xi$-map of $u$ at $x$} if
$\psi(x)=u(x)$ and for every cone $C_x$, there is a neighbourhood of $x$ in $\overline{\Om}$ such that, thereon, we have
 \beq \label{eq7.5a}
\big|\xi^\bot(u-\psi)\big|^2  \leq \,
C_x\big[-\xi^\top(u-\psi)\big].
 \eeq

\noi (2) \emph{The map $\psi \in C^2(\R^n)^N$ is a second contact $\xi$-map of $u$ at $x$} if
$\psi(x)=u(x)$ and for every cone $C_x$, there is a neighbourhood of $x$ in $\overline{\Om}$ such that, thereon, we have
 \beq \label{eq7.5b}
\big|\xi^\bot(u-\psi)\big|^2  \leq \,
(C_x)^2\big[-\xi^\top(u-\psi)\big].
 \eeq
\ed
In Definition \ref{de7.3} we allow for boundary points $x\in \p \Om$, but we are mostly interested in interior points $x\in \Om$. Inequalities \eqref{eq7.5a} and \eqref{eq7.5b} contain a lot of information. Specifically, \eqref{eq7.5b} says that for every $L>0$, exists an $r>0$ such that
 \begin{align}
\big|\xi^\bot(u-\psi)(y)\big|^2 \, \leq \,
L^2|y-x|^2\big[-\xi^\top(u-\psi)(y)\big] , \label{7.6b}
 \end{align}
for $y \in \mB_r(x)\cap \overline{\Om}$. Hence, \eqref{eq7.5b} is an elegant restatement of
 \begin{align}
\big|\xi^\bot(u-\psi)(y)\big|^2 \, \leq \,
o(|y-x|^2)\big[-\xi^\top(u-\psi)(y)\big] , \label{7.7b}
 \end{align}
as $\overline{\Om} \ni y \to x$, by using ``control by cones". Since the left hand side of \eqref{7.7b} is nonnegative and $u(x)=\psi(x)$, the projection $\xi^\top(u-\psi)$ along the line in $\R^N$ spanned by $\xi$ has a (local) \emph{vanishing maximum at $y=x$}:
\beq
\xi^\top(u-\psi)(y)\, \leq \, 0\, = \ \xi^\top(u-\psi)(x)
\eeq
for $y \in \overline{\Om}$ near $x$ (Figures 2(a), 3). 
\[
\underset{ \text{ Figure 2(a).\hspace{80pt} Figure 2(b). \hspace{80pt}Figure 2(c).} }{\includegraphics[scale=0.15]{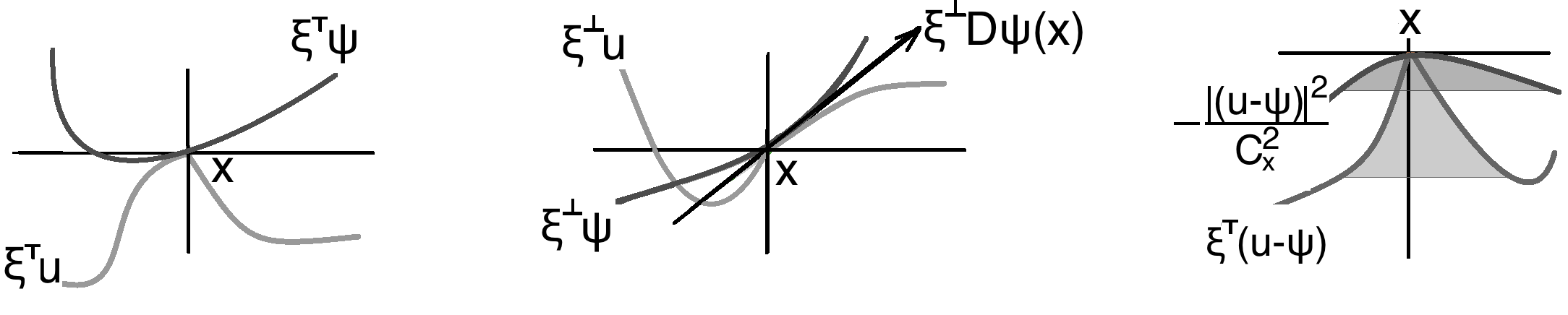}}
\]
Moreover, since the right hand side of \eqref{7.7b} is of order $o(|y-x|^2)$,  the projection $\xi^\bot(u-\psi)$ along the hyperplane $\xi^\bot \sub \R^N$ has a \emph{vanishing derivative at $y=x$}, since $\xi^\bot(u-\psi)(y)= o(|y-x|)$ as $y \to x$ (Figures 2(b), 3). Hence, since $u$ coincides with $\psi$ at $x$, the codimension-one projection $\xi^\bot u$ on the hyperplane is differentiable and $\D(\xi^\bot u)(x)= \xi^\bot \D\psi(x)$, although \textit{$\D u (x)$ may not exist}.
\[
\underset{\text{Figure 3.}}{\includegraphics[scale=0.28]{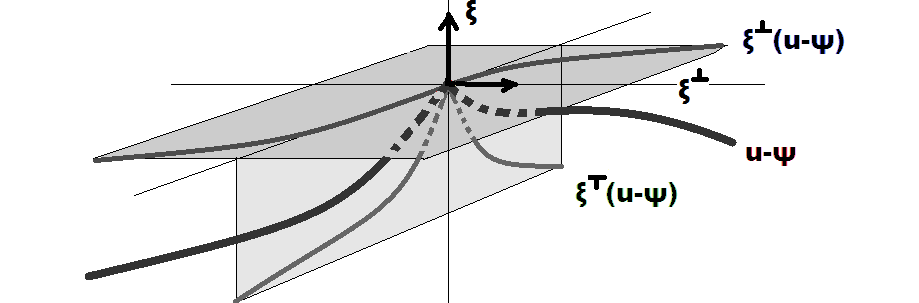}}
\]
Moreover, \eqref{7.7b} reveals that there is a coupling between $\xi^\top u$ and $\xi^\bot u$, which can be interpreted as that either the \emph{maximum of $\xi^\top(u-\psi)$ is constrained by $-|\xi^\bot(u-\psi)|^2/C_x^2$} or that the decay of \emph{$|\xi^\bot(u-\psi)|^2$ near $x$ is controlled by $-\xi^\top(u-\psi)$ via cones} (Figures 2(c), 3).

We will shortly see that contact maps  constitute an appropriate
notion of extremum for PDE theory. Let us first connect contact maps to contact jets. We will consider only the second order case and $x \in \Om$ and we refrain from providing details for the first order case and boundary points which can be done by simple modifications. Given a continuous $u : \R^n \supseteq \Om \larrow \R^N$, $x\in \Om$ and $\xi \in \mS^{N-1}$,  set
  \beq \label{eq7.6}
D^{2,\xi} u(x)\ := \ \left\{\big(\D\psi(x),\D^2 \psi(x)\big) \Bigg|
\begin{array}{l}
\psi \in C^2(\R^n)^N, \ \psi(x)=u(x)\ \ \& \ \, \forall \text{ cone } C_x,\\
\big|\xi^\bot(u-\psi)\big|^2 \! \leq (C_x)^2\big[ \!-\xi^\top(u-\psi)\big]
 \text{ near }x\\
\end{array}
\right\}.
  \eeq

 \bt[Equivalence between extremality and jets] \label{th7} If
 $u : \R^n \supseteq \Om \larrow \R^N$ is continuous, $x\in \Om$ and $\xi \in
 \mS^{N-1}$, then
  \beq \label{eq7.7}
D^{2,\xi} u(x) \, = \, J^{2,\xi} u(x).
  \eeq
That is, contact jets coincide with the set of derivatives of
contact maps.
 \et

The proof is based on the following lemma, which roughly states that
we can always absorb all of the second order Taylor remainder of a
contact $\xi$-map into its $\xi$-projection:

 \bl \label{le8} Let $\psi \in C^2(\R^n)^N$ be a second contact
 $\xi$-map of the continuous map $u : \R^n \supseteq \Om \larrow
 \R^N$ at $x\in \Om$. There exists a second contact $\xi$-map $\hat{\psi} \in C^2(\R^n)^N$ of $u$ at $x$ such that $\psi =\hat{\psi}$ up to second order
 at $x$ (that is, $\psi(x) =\hat{\psi}(x)$, $\D\psi(x)=D\hat{\psi}(x)$ and
  $\D^2\psi(x)=D^2\hat{\psi}(x)$) while the $\xi^\bot$-projection of
 the second order Taylor remainder of $\hat{\psi}$ vanishes.
 \el

\BPL \ref{le8}. Let $T_{2,x}$ and $R_{2,x}$ denote the operators of
second order Taylor polynomial and Taylor remainder at $x$
respectively. Then, by \eqref{eq7.5b}, we have
 \beq \label{eq7.8}
\Big|\xi^\bot\big(u-T_{2,x}\psi-R_{2,x}\psi \big)\Big|^2 \leq \,
C_x^2\Big[-\xi^\top\big(u-T_{2,x}\psi-R_{2,x}\psi \big)\Big]
 \eeq
locally in a neighbourhood of $x$. By employing Lemma
\ref{le3} for $N=1$, we can find an increasing $\rho\in
C^2(0,\infty)$ with $0\leq\rho(y)\leq o(|y-x|^2)$ such that $\rho \geq \xi^\top R_{2,x}\psi$ and
 \beq \label{eq7.9}
\Big|\xi^\bot\big(u-T_{2,x}\psi \big)-\xi^\bot R_{2,x}\psi\Big|^2
\leq \, \rho \Big[-\xi^\top\big(u-T_{2,x}\psi \big)\ +\ \rho \Big],
 \eeq
near $x$. By expanding the first term of \eqref{eq7.9}, we
estimate
 \begin{align} 
\Big|\xi^\bot\big(u-T_{2,x}\psi \big)\Big|^2\ &
\leq \ \rho \Big[-\xi^\top\big(u-T_{2,x}\psi \big)\ +\ \rho \Big]
 \, + \, \big|\xi^\bot R_{2,x}\psi\big|^2 
 \nonumber\\
& \ \ \ \, + \, \Big(\xi^\top\big(u-T_{2,x}\psi \big)\Big)^\top
\Big(\xi^\top\big( \xi^\bot R_{2,x}\psi \big)\Big) 
\nonumber
\\
& = \ \rho \Big[-\xi^\top\big(u-T_{2,x}\psi \big)\ +\ \rho \Big]
 \, + \, \big|\xi^\bot R_{2,x}\psi\big|^2 
\\
& \ \ \  \, + \, 2\Big(\frac{1}{\sqrt{2}}\xi^\top\big(u-T_{2,x}\psi
\big)\Big)^\top
\Big(\sqrt{2} \xi^\top\big( \xi^\bot R_{2,x}\psi \big)\Big). \nonumber
\\
&\leq \ \rho \Big[-\xi^\top\big(u-T_{2,x}\psi \big)\ +\ \rho \Big]
 \, + \, \big|\xi^\bot R_{2,x}\psi\big|^2 \\
& \ \ \  \, + \, \frac{1}{{2}} \Big|\xi^\top\big(u-T_{2,x}\psi
\big)\Big|^2 \, + \, 2\Big|\xi^\top\big( \xi^\bot R_{2,x}\psi
\big)\Big|^2. \nonumber
\\
 & \leq \ 2\rho
\Big[-\xi^\top\big(u-T_{2,x}\psi \big)\ +\ 2\rho \Big]
 \, + \, 2\big|\xi^\bot R_{2,x}\psi\big|^2. \nonumber
\end{align}
Hence,
 \begin{align} \label{eq7.11}
\Big|\xi^\bot\big(u-T_{2,x}\psi \big)\Big|^2\ 
 & = \ 2\rho
\Big[-\xi^\top\big(u-T_{2,x}\psi \big) \Big]\, + \, \Big\{4\rho^2
 \, + \, 2\big|\xi^\bot R_{2,x}\psi\big|^2\Big\} \nonumber\\
  & \leq  \ 2\big(\rho + |\xi^\bot R_{2,x}\psi| \big)
\Big[-\xi^\top\big(u-T_{2,x}\psi \big)  +  2\big(\rho + |\xi^\bot
R_{2,x}\psi| \big)\Big]. \nonumber
 \end{align}
In view of the above, the lemma follows by defining
$\hat{\psi} :=  T_{2,x}\psi \, + \, 2\big(\rho + |\xi^\bot
R_{2,x}\psi| \big) \xi$. Indeed, by construction we have $\psi(x) =\hat{\psi}(x)$, $\D\psi(x)=D\hat{\psi}(x)$, $\D^2\psi(x) =D^2\hat{\psi}(x)$ and $\xi^\bot
R_{2,x}\hat{\psi} \equiv 0$. Moreover, we have $\xi^\top R_{2,x}\hat{\psi} =  2\big(\rho +  |\xi^\bot R_{2,x}\psi| \big) \geq  \xi^\top R_{2,x}\psi$. Thus, we infer that
 \beq \label{eq7.13}
 \big|\xi^\bot(u-\hat{\psi})\big|^2 \, \leq\, \xi^\top R_{2,x}\hat{\psi}
 \big[-\xi^\top(u-\hat{\psi})\big],
 \eeq
and since $0\leq\xi^\top R_{2,x}\hat{\psi}(y)\leq o(|y-x|^2)$ as
$y\rightarrow x$, for every cone $C_x$ with vertex at $x$,
there exists a neighbourhood of $x$ such that, thereon,
 \beq \label{eq7.13}
 \big|\xi^\bot(u-\hat{\psi})\big|^2 \, \leq\, C_x^{2}
 \big[-\xi^\top(u-\hat{\psi})\big]. \ \ \  \qed 
 \eeq

We may now establish Theorem \ref{th7}.

\BPT \ref{th7}. Let $\psi$ be a second contact $\xi$-map of
$u$ at $x$. By Lemma \ref{le8}, there exists a second contact
$\xi$-map $\hat{\psi}$ of $u$ as $x$ such that
$\hat{\psi}=\psi$ up to second order at $x$ and moreover $\xi^\bot
R_{2,x}\hat{\psi} \equiv 0$. By \eqref{eq7.5b}, for any $L>0$, there
exists an $r>0$ such that
 \begin{align} \label{eq7.14}
& \Big|\xi^\bot\Big(u(z+x)-u(x)-
\D\psi(x)z-\frac{1}{2}\D^2\psi(x):z\ot z\Big)\Big|^2 
\nonumber\\
& \leq  L^2|z|^2 \Big[ \! -\xi^\top\Big(u(z+x)-u(x)- \D\psi(x)z-\frac{1}{2}\D^2\psi(x):z\ot z\Big) + \xi^\top R_{2,x}\hat{\psi}(z+x)\Big], \nonumber
 \end{align}
whenever $|z|\leq r$. By Lemma \ref{le4.3} for $N=1$, there exists
an increasing function $\si \in C^2(0,\infty)$ such that $\si(0^+)=0$ and $\si(|z|)|z|^2 \geq \xi^\top R_{2,x}\hat{\psi}$ as $z\rightarrow 0$, and also
\[
 \begin{split}
& \Big|\xi^\bot\Big(u(z+x)-u(x)- 
\D\psi(x)z-\frac{1}{2}\D^2\psi(x):z\ot z\Big)\Big|^2 
\\
& \leq \ \si(|z|)|z|^2 \Big[-\xi^\top\Big(u(z+x)-u(x)- \D\psi(x)z  -\frac{1}{2}\D^2\psi(x):z\ot z\Big)\ + \
\si(|z|)|z|^2\Big], 
 \end{split}
 \]
as $z\rightarrow 0$. By Theorem \ref{th6}, the above implies
$\big(\D\psi(x),\D^2\psi(x)\big) \in J^{2,\xi}u(x)$. Conversely,
let $(\P,\X) \in J^{2,\xi}u(x)$. Again by Theorem \ref{th6}, if
$\si$ is as stated, the map
 \beq \label{eq7.16}
\psi(z+x)\ :=\ u(x) \, + \, \P z \, +\, \frac{1}{2}\X:z\ot z\, +\, \si(|z|)|z|^2\xi
 \eeq
satisfies $\psi \in C^2(\R^n)^N$, $\psi(x)=u(x)$ and
 \beq \label{eq7.17}
 \big|\xi^\bot(u-{\psi})(y)\big|^2 \, \leq\, o(|y-x|^2)
 \big[-\xi^\top(u-{\psi})(y)\big],
 \eeq
as $y\rightarrow x$. Hence, $\psi$ is a second contact $\xi$-map of $u$ at $x$.
The theorem follows. \qed
\medskip

In view of Theorem \ref{th7}, we can reformulate Definitions
\ref{de2.3}-\ref{de2.2} of CS as follows (we state only the second order case for brevity):

 \bd[Contact solutions for second order systems, cf.\ Def.\ \ref{de2.3}]
 \label{de7.5}
\noi Suppose $\F$ is as in \eqref{3.1}. Then, the continuous map $u :\R^n \supseteq \Om \larrow \R^N$ is a \emph{contact solution to} \eqref{3.2} when for any $x\in \Om$ any $\xi \in \mS^{N-1}$ and any second contact $\xi$-map $\psi \in C^2(\R^n)^N$ of $u$ at $x$, we have
 \beq \label{eq7.20}
\xi^* \F\big(x,\psi(x),\D\psi(x),\D^2\psi(x)\big)\, \geq \, 0 .
  \eeq
\ed

The ``contact  principle calculus'' result we establish below explains why contact maps of the solution play the role of smooth ``test maps'' for the PDE system.

\bt[Nonlinear passage of derivatives to contact maps] \label{th8} Suppose $u : \R^n \supseteq \Om \larrow \R^N$ is a continuous map, twice differentiable at $x\in \Om$. Let $\psi \in C^2(\R^n)^N$
and $\xi \in \mS^{N-1}$ be given. Consider the following statements:

\smallskip

 \noi(i) $\psi$ is a second contact $\xi$-map of $u$ at $x \in \Om$.

 \noi (ii) We have
\[
\left\{ \ \ \
\begin{split}
\D\big( u-\psi \big)(x) & = \ 0, \\
\xi \vee \D^2\big( u- \psi \big)(x) & \leq_\ot  \!0. 
\end{split}
\right.
\]
Then, (i) implies (ii). Moreover, (ii) implies (i) if moreover $\xi^\top \D^2\big(u-\psi \big)(x) < 0$.
 \et
Trivial modifications in the arguments of the proof of Theorem
\ref{th8} that follows readily imply the following consequence.

 \begin{corollary}[first order contact] \label{cor6}  In the setting of Theorem \ref{th8},
  we have that if $u$ is differentiable, then $\psi \in C^1(\R^n)^N$ is a first contact $\xi$-map
of $u$ at $x$ if and only if $\D\big(u-\psi \big)(x)=0$.
 \end{corollary}

Theorem \ref{th8} has already been established implicitly, in the
language of contact jets. Indeed, one may employ Theorem \ref{th7} and Proposition
\ref{th3} to remove the disguise. Further, one can also easily establish the next consequence.

 \begin{corollary}[Rank-One Decompositions] \label{cor7}
  In the setting of Theorem \ref{th8},
we have that if $\psi$ is a second contact $\xi$-map of $u$
at $x$, we have the rank-one decompositions
 \begin{align}
\D(u-\psi)(x)  & = \xi \ot  \D \big(\xi^\top(u  -\psi)\big)(x) ,   \label{eq7.36} \\
\D^2(u-\psi)(x) & =  \xi \ot  \D^2\big(\xi^\top (u -\psi)\big)(x) ,
\label{eq7.37}
\\
  \D \big(\xi^\top(u & -\psi)\big)(x)  = 0, \label{eq7.38}\\
 \D^2\big(\xi^\top (u &-\psi)\big)(x)  \leq  0. \label{eq7.39}
 \end{align}
Conversely, if the relations \eqref{eq7.36}-\eqref{eq7.39} hold true, then $\psi$ is
a second contact $\xi$-map of $u$ at $x$ if in addition the inequality \eqref{eq7.39} is strict.
 \end{corollary}

\section{ Approximation and stability of contact jets} 
\label{section8}

In this section we consider the problem of stability of CS under limits. We recall that in the scalar case, VS pass to limits under merely locally uniform convergence. This means that if the sequence of solutions $(u_j)_1^\infty$ to equations $F_j(\cdot,u_j,\D u_j,\D^2 u_j)=0$ satisfies $u_j \to u$ and the nonlinearities satisfy $F_j \to F$, both convergences locally uniform as $j \to \infty$, then $u$ solves $\F(\cdot,u,\D u ,\D^2 u   )=0$. This important property is a consequence of the fact that maxima perturb to maxima under uniform convergence. 

We begin with a counterexample which shows that in the vectorial case this property \emph{fails}. More precisely, \emph{second contact maps do not perturb to second contact maps not even under strong $C^1$ convergence}, and neither $C^{1,\al}$ convergence suffices for any $\al <1$. Moreover, the first order variant of Example \ref{ex49} below shows that \emph{first contact maps do not perturb to first contact maps under $C^\al$ convergence} for any $\al <1$.

\begin{example}[Instability of contact maps] \label{ex49} For any $\al \in (0,1)$, there exist $u\in C^{1,\al}(\R)^2$, $\xi \in \mS^1$, a second contact $\xi$-map $\psi \in C^2(\R)^2$ of $u$ at $x=0$ and a sequence $(u_m)^\infty_1 \sub C^\infty(\R)^2$ such that $\xi^\top u_m \to \xi^\top u$ in $C^1(\R)$ and $\xi^\bot u_m \to \xi^\bot u$ in $C^\infty(\R)^2$ as $m\to \infty$, but there exists no sequence of second contact $\xi$-maps $\psi_m \in C^2(\R)^2$ of $u_m$ along any $x_m\to 0$ such that
$\big(\D\psi_m(x_m), \D^2\psi_m(x_m) \big) \to \big(\D\psi(0), \D^2\psi(0)\big)$ as $m\to \infty$.

Indeed, define $u$ by $u(z):=(-|z|^{1+\al},0)^\top$, fix $k\in \R\set \{0\}$ and set $\xi:=(1,0)^\top=e_1$ and $\psi(z):=(0,k|z|^2)^\top$. We first verify that $\psi$ is a second contact $e_1$-map at $x=0$: for any $L>0$, there is $r>0$ such that for $|z|< r$,
\begin{align}
\big|\xi^\bot(u-\psi)(z) \big|^2 = \big|0-k|z|^2 \big|^2 =  \big(k|z|^{\frac{1-\al}{2}}\big)^2|z|^2 \left[-(-|z|^{1+\al} -0)\right]  \nonumber
\end{align}
which yields $\big|\xi^\bot(u-\psi)(z) \big|^2\leq\ L^2|z|^2\big[-\xi^\top(u-\psi)(z)\big]$. Consider now the sequence $u_m:= \eta^{1/m}*u$ where $\eta^{1/m}$ is the standard mollification of $u$. Since $\xi^\bot u\equiv 0$, we obtain that $u_m= (\eta^{1/m}*(\xi^\top u), 0 )^\top$ and consequently (a) and (b) follow. Since $u_m$ is smooth, it possesses contact maps at all $x\in \R$. Choose a sequence $x_m \to 0$ and let $\psi_m$ be a second contact $\xi$-map of $u_m$ at $x_m$. Then, by Corollary \ref{cor3} we have that $\xi^\bot \psi_m$ is equal up to second order to $\xi^\bot u_m$ at $x_m$. Hence, as $m\to \infty$ we have
\beq
D^2 \xi^\bot \psi_m (x_m)\, = \, D^2\xi^\bot u_m (x_m)\, =\, 0\ \ \not\!\!\longrightarrow \ 2ke_2 \, = \, D^2\xi^\bot \psi(0).
\eeq
\end{example}
The reason of instability of contact maps is that projections on the normal hyperplane $\xi^\bot$ of ``generalised hessians" can be the whole tensor subspace $\xi^\bot \ot \R^{n\by n}_s \sub \R^N \! \ot \R^{n\by n}_s $. This phenomenon is a general fact, which appears when $\xi^\top u$ fails to be $C^{1,1}$ near the basepoint. The next lemma exhibits the previous situation and supplements Corollary \ref{cor3}. By utilising this result and the properties of the Little H\"older space $c^{1,\al}(\R^n)^N$ which is the closure of $C^\infty_c(\R^n)^N$ under the H\"older norm, it follows that not even $C^{1,\al}$ convergence suffices (if $\al<1$).

\bl \label{le50} Let $\psi \in C^2(\R^n)^N$ be a second contact $\xi$-map of the continuous map $u : \R^n \supseteq \Om \larrow \R^N$ at $x\in \Om$ for some $\xi \in \mS^{N-1}$. We set:
\[
l^-:=\, \liminf_{y\to x} \frac{-\xi^\top(u-\psi)(y)}{|y-x|^2}\ , \ \ \ l^+:=\, \limsup_{y\to x} \frac{-\xi^\top(u-\psi)(y)}{|y-x|^2}.
\] 
If $l^-=\infty$, then, all quadratic perturbations $\hat{\psi}(y):=\psi(y) + \frac{1}{2}\xi^\bot \X: (y-x) \ot (y-x)$ are also contact $\xi$-maps for any $\X \in \R^N \! \ot\R^{n\by n}_s$.
If $l^+<\infty$, then $\xi^\bot \psi$ is unique (up to superquadratic perturbations $o(|y-x|^2)$ as $y \to x$).
\el

\BPL \ref{le50}. If $l^-=\infty$, there is $\om \in C^0(0,\infty)$ with $\om > \om(0^+)=0$ with
\[
\frac{-\xi^\top(u-\psi)(y)}{|y-x|^2}  \geq \frac{1}{\om^2(|y-x|)}
\]
for $|y-x|<1$. Hence, $-\xi^\top(u-\psi)(y)\, \om^2(|y-x|)|y-x|^2 \geq |y-x|^4$ and since $\psi$ is a contact $\xi$-map, we estimate as $y\to x$
\begin{align}
\big|\xi^\bot(u-\hat{\psi})(y)\big|^2\ &\leq \ 2 \big|\xi^\bot(u- {\psi})(y)\big|^2\ +\  \frac{1}{2}\big|\xi^\bot \X: (y-x) \ot (y-x) \big|^2 \nonumber\\
&\leq \ \Big( o(1) \, +\, \frac{1}{2} \big\|\xi^\bot \X \big\|^2 \om^2(|y-x|) \Big)|y-x|^2  \big[-\xi^\top(u-\psi)(y)\big] . \nonumber
\end{align}
To see the last claim, apply Corollary \ref{cor3} and  Lemma \ref{le8}.    \qed 
\ms

Happily enough, the discouraging instability of contact maps is not detrimental to the stability of CS. The reason is that when we try to approximate a system by adding a ``viscosity term", there exists some extra information which is trivial in the scalar case and allows convergence of the approximating solutions. In order to make this statement precise, we introduce an auxiliary notion of sequential derivatives needed in the exploitation of stability and approximation.

\bd[Approximate derivative] \label{de51} let $u:\R^n \supseteq \Om \larrow \R^N$ be a continuous map. The set of \emph{Approximate first jets} of $u$ at $x\in \Om$ is
\beq \label{9.8}
A^1u(x)\ :=\ \Big\{\P \in  \R^N \! \ot\R^n \ \Big| \ \liminf_{r\to 0}\max_{|z|=r}\frac{\big|u(z+x)-u(x)-\P z \big|}{r} =0 \Big\}
\eeq
The set of \emph{Approximate second derivatives} of $u$ at $x\in \Om$ is
\begin{align} \label{9.9}
A^2u(x)\ :=\ \Big\{ & (\P,\X)  \in \ \R^N \! \ot\big(\R^n \by \R^{n\by n}_s\big)\ \Big| \nonumber\\
& \liminf_{r\to 0}\max_{|z|=r}\frac{\big|u(z+x)-u(x)-\P z-\frac{1}{2}\X:z \ot z\big|}{r^2}=0 \Big\}
\end{align}
\ed
\begin{remark} 
Obviously, if $u$ is (twice) differentiable at $x$, then $A^1u(x)=\{\D u (x)\}$ and $A^2u(x)=\{(\D u (x),\D^2 u(x))\}$. In general, approximate derivatives may exist at non-differentiability points, as it happens for the Lipschitz continuous function $u :\R \to \R$ given by $u(z):=z\cos({1}/{|z|})$ for $z\neq 0$ and $u(0)=0$ for which $[-1,+1]= A^1u(0)\neq \emptyset$, while $u'(0)$ does not exist. This follows from the observation $\max_{|z|=r}(|u(z)-u(0)- pz|/r)=|\cos(1/r) -p|$.
\end{remark}
The following is the main approximation result for contact jets. It follows that \emph{contact jets perturb to contact jets under weak$^*$ convergence in the local Lipschitz space}, together with a \emph{technical assumption} which appears to be satisfied in the cases of interest. This assumption requires convergence of codimension-one projections of sequential jets along a sequence of hyperplanes.

\bt[Approximation of contact jets] \label{th54} Let $u : \R^n \supseteq \Om \larrow \R^N$ be continuous and fix $(\P,\X) \in J^{2,\xi}u(x)$ for some $x\in \Om$, $\xi \in \mS^{N-1}$. Suppose there exists $(u_m)_1^\infty \sub C^2(\Om)^N$ such that $u_m \weakstar  u $ in $W^{1,\infty}_{loc}(\Om)^N$ and for some $e \not\in \xi^\bot$ we have
\beq \label{9.11}
\left\{
\begin{array}{l}
\text{For any $y_m \to x$, there is $(e_m)^\infty_1 \sub \mS^{N-1}$ such that}\\
\text{$e_m \to e$  \ and }\  \dist\Big( e_m^\bot \D u _m(y_m)\, , \, A^1(e^\bot u)(x)\Big) \to  0.
\end{array}
\right.
\eeq
Then, there exist sequences $(x_m)_1^\infty$ and $(\P_m,\X_m)^\infty$ with $(\P_m,\X_m) \in J^{2,\xi}u_m(x_m)$ satisfying $x_m \to x$ and also $\big(\P_m , \xi^\top \X_m\big) \larrow \big(\P , \xi^\top \X\big)$ as $m \to \infty$.
\et
Theorem \ref{th54} is optimal: by example \ref{ex49}, not even $C^1$ convergence $u_m \to u$ suffices to guarantee $\xi^\bot \X_m \to \xi^\bot \X$. There is a ``loss of information" which occurs when $e \in \xi^\bot$ (i.e.\ when $e^\top \xi=0$). The proof is based on the next result which relates approximate derivatives of codimension-one projections to contact jets.

\begin{proposition}[Approximate derivatives and contact jets on hyperplanes] \label{Pr55} Let $u : \R^n \supseteq \Om \larrow \R^N$ be continuous and $x\in \Om$ and $e,\, \xi \in \mS^{N-1}$. Then, we have:
\beq
\left. \label{9.13}
\begin{array}{l}
\P\, \in\,  J^{1,\xi}u(x) \\
\Q\, \in \,A^1(e^\bot u)(x)
\end{array}
\right\}
\ \ \Longrightarrow\ \ \Q = e^\bot \P \, , \  \text{ if\  } e\neq \pm \xi.
\eeq
In particular, if both sets $e^\bot \big( J^{1,\xi}u(x)\big)$, $A^1(e^\bot u)(x)$ are nonempty, they are singletons and coincide. If moreover $J^{2,\xi}u(x)\neq \emptyset$, then $\Q = e^\bot \P$. Further:
\beq
\left.  \label{9.14}
\begin{array}{l}
(\P,\X)\, \in\,  J^{2,\xi}u(x) \\
(\Q,\Y)\, \in \,A^2(e^\bot u)(x)
\end{array}
\right\}
\ \ \Longrightarrow\ \ \left\{
\begin{array}{r}
\Q - e^\bot \P\, = \ \, 0, \\
(e^\bot \xi) \vee [\Y - e^\bot \X]\leq_\ot \! 0.
\end{array}
\right.
\eeq
\end{proposition}

\begin{corollary}
By \eqref{9.14} and Lemma \ref{le2}, we deduce that $\Y=e^\bot \X$ on the hyperplane $(e^\bot \xi)^\bot$ of $\R^N$ and $(e^\bot \xi)^\top (\Y - e^\bot \X)\leq 0$ along $(e^\bot \xi)\ot (e^\bot \xi)$ (Figures 6(a),(b)).
\[
\underset{\text{Figure 4(a): Illustration for $N=2$\hspace{60pt} Figure 4(b):  Illustration for $N=3$\ \ \ \ \ \ \ }}{\includegraphics[scale=0.4]{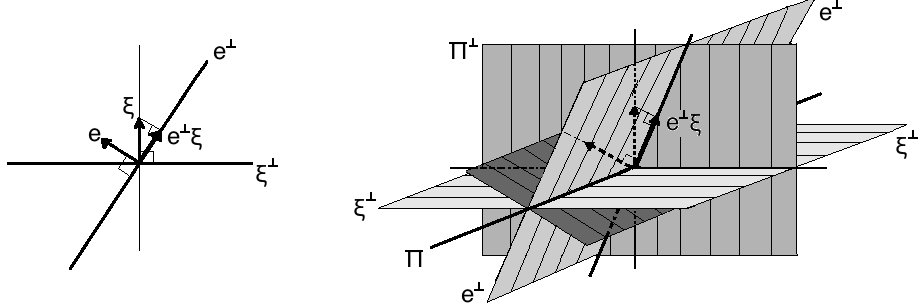}}
\]
\end{corollary}

\begin{remark} 
Proposition \ref{Pr55} is optimal, as the example $u(z):=-|z|\cos^2(1/z)\chi_{\R\set\{0\}}$ shows: indeed, we have $J^{1,+}u(0)=\{0\}=A^1u(0)$ but $u'(0)$ does not exist, although $J^{1,+}u(0)$ and $A^1u(0)$ coincide and are singletons. Namely, some ``loss of information" occurs when $e=\pm \xi$ (i.e.\ when $e^\bot \xi =0$).
\end{remark}

\BPP \ref{Pr55}. Since $\Q \in A^1(e^\bot u)(x)$, exists $r_j \to 0$ such that
\beq \label{9.15}
e^\bot \Big(u(wr_j +x)-u(x)-r_j \Q w\Big)\, = \, o(r_j)
\eeq
as $j \to \infty$, for all $w\in \mS^{n-1}$. Fix $\theta \in e^\bot$. If $e\neq \pm \xi$, then $e^\bot \set \xi^\bot \neq \emptyset$, and for any $\e>0$, exists $\theta_\e \in e^\bot \set \xi^\bot$ with $|\theta - \theta_\e|\leq \e$. Since $(\theta_\e \ot \theta_\e)e^\bot = \theta_\e \ot \theta_\e$,  \eqref{9.15} gives
\begin{align} \label{9.16}
o(r_j) = -\, (\theta_\e^\top \xi) \left[\theta_\e^\top  \Big(u(wr_j +x)-u(x)-r_j \Q w\Big)\right]
\end{align}
as $j \to \infty$. Since $\P \in J^{1,\xi}u(x)$, we have
\begin{align} \label{9.17}
o(|z|) =  (\theta_\e^\top \xi) \left[\theta_\e^\top  \Big(u(z +x)-u(x)-(e^\bot \P) z\Big)\right]
\end{align}
as $z\to 0$. By choosing $z=w r_j $ in \eqref{9.17} and summing \eqref{9.17} and \eqref{9.16}, we get $(\theta_\e^\top \xi)\big[ \theta_\e^\top \big(\Q-e^\bot \P\big)w\big] \leq o(1)$ as $r_j\to 0$. By interchanging $w$ with $-w$, using that $|\theta_\e^\top \xi|>0$ and letting $j\to \infty$, we get  $\theta_\e^\top \big(\Q-e^\bot \P\big)=0$. By letting $\e \to 0$, we find $\theta^\top \big(\Q-e^\bot \P\big)=0$. Since $\theta \in e^\bot$ is arbitrary and $\Q = e^\bot\Q$, we conclude that $\Q = e^\bot\P$. If $e=\pm \xi$ and moreover $J^{2,\xi}u(x)\neq \emptyset$, Corollary \ref{cor3} and Definition \ref{de51} imply that $\Q = e^\bot\P=\D(e^\bot u)(x)$. Further, if $(\Q,\Y) \in A^2(e^\bot u)(x)$, it trivially follows that $\Q \in A^1(e^\bot u)(x)$. Since $(\P,\X) \in J^{2,\xi}(x)\neq \emptyset$, part (a) implies that $\Q = e^\bot\P$. Fix $\theta \in e^\bot$. By arguing as eaelier, there exists $r_j \to 0$ such that
\begin{align} 
o(r_j^2)\ &\geq \ \theta \ot \theta : \Big[\xi \vee  \Big(u(wr_j +x)-u(x)- r_j\Q w -\frac{r_j^2}{2} (e^\bot\X): w\ot w\Big)\Big], \label{9.19}\\
o (r_j^2)\ &\geq  \ - \theta \ot \theta : \Big[\xi \vee  \Big(u(wr_j +x)-u(x) -r_j \Q w- \frac{r_j^2}{2} \Y: w\ot w\Big)\Big], \label{9.20}
\end{align}
as $j\to \infty$, for all $w\in \mS^{N-1}$. By writing $\theta = e^\bot \eta$ for some $\eta \in \R^N$, using the symmetry of $e^\bot$ and summing \eqref{9.19} and  \eqref{9.20}, we obtain as $r_j \to 0$
\begin{align}
o(1) \ &\geq \ (e^\bot \eta) \ot (e^\bot \eta) : \Big[\xi \vee \big(\Y-e^\bot \X\big): w\ot w\Big] \nonumber\\
 &=  \Big[(e^\bot \xi) \vee (\Y-e^\bot \X )\Big]: \eta \ot w \ot \eta \ot w.\nonumber
\end{align}
By passing to the limit we conclude that $(e^\bot \xi) \vee (\Y-e^\bot \X )\leq_\ot 0$.     \qed

\ms

\BPT \ref{th54}. Since $(\P,\X) \in J^{2,\xi}u(x)$, we have $(\xi^\top \P, \xi^\top\X) \in J^{2,+}(\xi^\top u)$ $(x)$. By the $C^0$ convergence $\xi^\top u_m \to \xi^\top u$ as $m\to \infty$, standard arguments of the scalar case (see e.g.\ \cite{CIL,K}) imply that there exists $x_m \to x$ and $(p_m,X_m) \in J^{2,+}(\xi^\top u_m)(x_m)$ such that 
\beq \label{9.22}
(p_m,X_m) \larrow (\xi^\top \P,\xi^\top \X)\ ,\ \ \text{ as }m\to \infty.
\eeq
Since $u_m \in C^2(\Om)^N$, it follows that $p_m = \xi^\top \D u_m(x_m)$ and $X_m \geq \xi^\top \D^2 u _m(x_m)$. By Theorem \ref{th2}, the set $ J^{2,\xi}u_m(x_m)$ contains
\beq   \label{9.24}
(\P_m,\X_m)\, :=\, \Big(\D u _m(x_m),\,  \D^2(\xi^\bot u_m)(x_m)+ \xi \ot X_m\Big).
\eeq
By decomposing $\P_m =\xi \ot p_m +\xi^\bot \D u _m(x_m)$, in view of \eqref{9.22} and  \eqref{9.24} we see that it suffices to use \eqref{9.11} in order to show that $\xi^\bot \P_m \to \xi^\bot \P$ as $m\to \infty$.  For the sequence $x_m \to x$, assumption \eqref{9.11} implies that there exists a convergent sequence $\mS^{N-1} \ni e_m \to e$ of directions and an Approximate jet $\Q \in A^1(e^\bot u)(x)$ such that $e_m^\bot \D u _m(x_m) \to \Q$ as $m\to \infty$. Since $(\P,\X) \in J^{2,\xi}u(x)$, Proposition \ref{Pr55} implies that $\Q =e^\bot \P$ and as a result we deduce that $e_m^\bot \P_m \to e^\bot \P$  as $m\to \infty$. Further, by replacing as we can  $\{(e_m)_1^\infty,e\}$ by $\{(-e_m)_1^\infty,-e\}$, we may assume that $0\leq \xi^\top e \leq 1$, namely that $\xi,e$ lie in the same halfspace (Figures 6(a),(b)). We distinguish two cases:

\textbf{Case 1:  $0<\xi^\top e <1$}. We use expansions with respect to non-orthonormal coordinates in order to show that $\xi^\bot \P_m \to \xi^\bot \P$  as $m\to \infty$. We define the codimension-two subspaces $\Pi_m := \xi^\bot \cap e_m^\bot$ and $\Pi := \xi^\bot \cap e^\bot$ which are intersections of hyperplanes (see Figures 6(a),(b)) and allow to write 
\beq \label{9.27}
\R^N \ =\  \Pi_m \oplus \spn[\{\xi, e_m^\bot \xi\} ]\ =\ \Pi \oplus \spn[\{\xi, e^\bot \xi\} ].
\eeq
Let us now de define the unit vectors $\eta_m := \sgn( e_m^\bot \xi)$ and $\eta := \sgn( e^\bot \xi)$. By expansion on the non-orthonormal frames $\{\xi, \eta_m, \Pi_m\}$, $\{\xi, \eta, \Pi\}$, we have that for any $a\in \R^N$, there exists $ \la_m(a), \mu_m(a),\la(a), \mu(a) \in \R$ such that
\begin{align}
a\ &= \ \la_m(a)\xi \ +\ \mu_m(a)\eta_m\, + \, \Pi_m a, \label{9.29}\\
a\ &= \ \la(a)\xi \ +\ \mu(a)\eta\ +\ \Pi a. \label{9.30}
\end{align}
Since $\xi,\eta_m$ are normal to $\Pi_m$ and $|\xi|=|\eta_m|=1$, by projecting \eqref{9.29}-\eqref{9.30} along $\xi \ot \xi $ and $\eta_m \ot \eta_m $ we obtain
\beq
\left\{
\begin{array}{l} \label{9.31}
\eta_m^\top a\, = \, \la_m(a)(\xi^\top \eta_m) \ +\ \mu_m(a), \ms\\
\xi^\top a\, = \, \la(a) \, +\, \mu_m(a) (\xi^\top \eta_m). 
\end{array}
\right.
\eeq
By solving the linear system \eqref{9.31}, we find
\beq
\begin{array}{l} \label{9.32}
\la_m(a)=    \dfrac{ \xi^\top a\, - \,  (\xi^\top \eta_m) (\eta_m^\top a)}{1\, -\, (\xi^\top \eta_m)^2}, \ \ \ 
\mu_m(a) =  \dfrac{ \eta_m^\top a\, - \,  (\xi^\top \eta_m) (\xi^\top a)}{ 1\, -\, (\xi^\top \eta_m)^2}, . 
\end{array}
\eeq
We observe that $\xi^\top (\xi^\bot \P_m)=0$ and use \eqref{9.29} and \eqref{9.32} to expand  $\xi^\bot \P_m$ as 
\beq \label{9.33}
\xi^\bot \P_m\, = \, \left[\frac{ -  (\xi^\top \eta_m) \xi \ot\eta_m}{ 1 - (\xi^\top \eta_m)^2 }\, +\, \frac{ \eta_m \ot\eta_m }{ 1-(\xi^\top \eta_m)^2 }\, +\, \Pi_m \right]\xi^\bot \P_m.
\eeq
By recalling that $\Pi_m = e_m^\bot \cap \xi^\bot$, we observe that $\Pi_m \xi^\bot =\Pi_m$. By using this and that $\xi^\bot =I-\xi \ot \xi$, we rewrite \eqref{9.33} as
\beq \label{9.35}
\begin{split}
\xi^\bot \P_m\, & = \, \left[\frac{ -  (\xi^\top \eta_m) \xi  + \eta_m}{ 1 - (\xi^\top \eta_m)^2 } \right] \ot \Big[\eta_m^\top \P_m -  (\xi^\top \eta_m)(\xi^\top \P_m)   \Big] \, +\, \Pi_m\P_m 
\\
& = \, \frac{\xi^\bot \eta_m}{ |\xi^\bot \eta_m|^2 } \ot \Big[\eta_m^\top \P_m -  (\xi^\top \eta_m)(\xi^\top \P_m)   \Big] \, +\, \Pi_m\P_m. 
\end{split}
\eeq
Similarly, we have
\beq \label{9.36}
\xi^\bot \P\, = \, \frac{\xi^\bot \eta}{ |\xi^\bot \eta|^2 } \ot \Big[\eta^\top \P -  (\xi^\top \eta)(\xi^\top \P)   \Big] \, +\, \Pi \P. 
\eeq
By using that $e_m \to e$, we obtain $e_m^\bot = I-e_m \ot e_m \larrow I-e \ot e = e^\bot$, as $m\to \infty$. Also, since $e_m^\bot \xi \to e^\bot \xi$ and $|e_m^\bot \xi|\geq |e^\bot \xi|/2 >0$ for $m$ large enough, we have $\eta_m \larrow \eta$. Since by assumption $e\not\in \xi^\bot$, for $m$ large we have $|\xi^\bot \eta_m| \geq |\xi^\bot \eta|/2>0$. Hence, $\sgn(\xi^\bot \eta_m) \larrow \sgn(\xi^\bot \eta)$ as $m\to \infty$. By \eqref{9.27} we have $e_m^\bot =\Pi_m +\eta_m \ot \eta_m $ and hence
$\Pi_m = e_m^\bot - \eta_m \ot \eta_m\ \larrow \ e_m^\bot - \eta \ot \eta= \Pi$ as $m\to \infty$. Since $\xi^\top \P_m \to \xi^\top \P$ and $\xi^\top \eta_m \to \xi^\top \eta$, in view of \eqref{9.35}, \eqref{9.36}, it suffices to show that $\Pi_m \P_m \to \Pi \P$ and that  $\eta_m^\top \P_m \to \eta^\top \P$ as  $m\to \infty$. Indeed we have the estimate
\begin{align}
\big| \Pi_m \P_m - \Pi \P\big|
&\, \leq\, \big| \Pi_m (e_m^\bot \P_m -e^\bot \P) \big|\, +\,  \big| (e^\bot \P)^\top (\Pi_m -\Pi)\big|\ \larrow \ 0,  \nonumber
\end{align}
 as $m\to \infty$, and similarly we conclude that $\eta_m^\top \P_m \to \eta^\top \P$.

\ms
\textbf{Case 2:  $\xi = e$}. If in addition $\xi =e=e_m$ for infinitely many terms $(m_k)_1^\infty$, then $\xi^\bot \P_{m_k} =e_{m_k}^\bot \P_{m_k} \to e^\bot \P = \xi^\bot \P$ as $k\to \infty$ and the   conclusion follows. If on the other hand $\xi^\top e_m <1$ for $m$ large enough, then the arguments of the previous case fail because we see that $\eta_m \not\to \eta$. Instead, by using that $\P_m =\D u _m(x_m)$ and $x_m \to x$, the local weak$^*$ convergence $u_m \weakstar u$ gives for $m$ large enough the bound
\beq \label{9.39}
|\P_m|\, \leq\, \big\| \D u _m\big\|_{L^\infty (\mB_R(x))}  \leq \, C(R)\ , \ \ \ R\, :=\, \frac{1}{2}\dist(x,\p \Om).
\eeq
By using \eqref{9.39} and that $\xi = e$, we estimate
\begin{align}  \label{9.40}
\big| \xi^\bot \P_m - \xi^\bot\P \big|\
& = \ \Big| \xi^\bot \big( e_m \ot (e_m^\top \P_m)\, +\, e_m^\bot \P_m \big)\ -\ \xi^\bot\P \Big| \nonumber\\
& \leq \ |\xi^\bot e_m| \, \big|e_m^\top \P_m \big|\  +\ \big|\xi^\bot(e_m^\bot \P_m -  \xi^\bot\P )\big|\\
& \leq \ |e^\bot e_m| \Big( |e_m\!-e| \, |\P_m| + \big|\xi^\top \P_m \big|\Big)\  +\ \big| e_m^\bot \P_m -  e^\bot \P \big|.
\nonumber
\end{align}
Since $e_m \to e$, $e_m^\bot e \to 0$, $e_m^\bot \P_m \to e^\bot \P$ and $\xi^\top \P_m \to \xi^\top \P$, the bounds \eqref{9.39} and \eqref{9.40} allow us to conclude.           \qed

\bibliographystyle{amsplain}

\end{document}